\input amstex
\documentstyle{amsppt}
\magnification=\magstep1
\vsize =21 true cm
\hsize =16 true cm
\loadmsbm
\topmatter

\centerline{\bf Icosahedron, exceptional singularities and modular forms}
\author{\smc Lei Yang}\endauthor
\endtopmatter
\document

\centerline{\bf Abstract}

\vskip 0.5 cm

  We find that the equation of $E_8$-singularity possesses two distinct
symmetry groups and modular parametrizations. One is the classical icosahedral
equation with icosahedral symmetry, the associated modular forms are theta constants
of order five. The other is given by the group $\text{PSL}(2, 13)$, the associated
modular forms are theta constants of order $13$. As a consequence, we show that
$E_8$ is not uniquely determined by the icosahedron. This solves a problem of Brieskorn
in his ICM 1970 talk on the mysterious relation between exotic spheres, the icosahedron
and $E_8$. Simultaneously, it gives a counterexample to Arnold's $A, D, E$ problem, and
this also solves the other related problem on the relation between simple Lie algebras
and Platonic solids. Moreover, we give modular parametrizations for the exceptional
singularities $Q_{18}$, $E_{20}$ and $x^7+x^2 y^3+z^2=0$ by theta constants of order $13$,
the second singularity provides a new analytic construction of solutions for the Fermat-Catalan
conjecture and gives an answer to a problem dating back to the works of Klein.

\vskip 0.5 cm

\centerline{\bf Contents}
$$\aligned
 &\text{1. Introduction}\\
 &\text{2. Standard structure on the $E_8$-singularity: the icosahedron}\\
 &\text{3. Exotic structure on the $E_8$-singularity: $\text{PSL}(2, 13)$}\\
 &\text{4. Modular parametrizations for some exceptional singularities}\\
 &\text{\quad and Fermat-Catalan conjecture}
\endaligned$$

\vskip 0.5 cm

\centerline{\bf 1. Introduction}

\vskip 0.5 cm

  The icosahedron, and the other regular solids, have been known as the Platonic solids since
the time of the ancient Greeks and have played an important role in the development of mathematics.
Moreover they turn out to have unexpected relations to many other topics in mathematics. Some of the
beauty and fascination associated with the icosahedron are given by the following correspondence
between the finite subgroups of $\text{SL}(2, {\Bbb C})$ and the $ADE$ type Dynkin diagrams (see \cite{Sl2}):
$$\matrix
  \text{Group} & \text{Root system}\\
  \text{the cyclic group of order $n$} & A_{n-1}\\
  \text{the binary dihedral group of order $4n$} & D_{n+2}\\
  \text{the binary tetrahedral group of order $24$} & E_6\\
  \text{the binary octahedral group of order $48$} & E_7\\
  \text{the binary icosahedral group of order $120$} & E_8
\endmatrix\eqno{(1.1)}$$
This correspondence appears in many branches of mathematics and physics, such as singularity theory (see
\cite{Ar1}, \cite{Ar3}, \cite{Br3}, \cite{Sl1}, \cite{Sl2}, \cite{Dur}), modular forms and arithmetic (see \cite{Hi3},
\cite{Na}, \cite{Se}, \cite{Du}), algebraic geometry (see \cite{B}, \cite{DV}, \cite{GV}, \cite{IN}, \cite{R})), algebra
(see \cite{DM}), representation theory (see \cite{Ko1}, \cite{Ko2}, \cite{Mc1}, \cite{Mc2}, \cite{N}), Lie
groups (see \cite{V}), geometry (see \cite{At}, \cite{H}), geometric topology (see \cite{KS}), differential
geometry (see \cite{Kr}), conformal field theory and subfactors (see \cite{CIZ}, \cite{J1}, \cite{J2}, \cite{J3},
\cite{Z}) and string theory (see \cite{W}). A connection between the left column and the right column in (1.1)
is given by the singularities (see \cite{Mc1}). In fact, the relation between simple singularities and simple
Lie groups is one of the most beautiful discoveries in mathematics and it will lead back to the classification
of Platonic solids (see \cite{B}). Starting with the polynomial invariants of the finite subgroup of
$\text{SL}(2, {\Bbb C})$, a surface is defined from the single syzygy which relates the three polynomials in two
variables. This surface has a singularity at the origin; the singularity can be resolved by constructing a smooth
surface which is isomorphic to the original one except for a set of component curves which form the pre-image
of the origin. The components form a Dynkin curve and the matrix of their intersections is the negative
of the Cartan matrix for the appropriate Lie algebra. The Dynkin curve is the dual of the Dynkin graph.
For example, if $\Gamma$ is the binary icosahedral group, the corresponding Dynkin curve is that of $E_8$,
and ${\Bbb C}^2/\Gamma \subset {\Bbb C}^3$ is the set of zeros of the equation $x^2+y^3+z^5=0$ which can
be parametrized by theta constants of order five (see section two for more details). The link of this
$E_8$-singularity, the Poincar\'{e} homology $3$-sphere (see \cite{KS}), has a higher dimensional lifting:
$z_1^5+z_2^3+z_3^2+z_4^2+z_5^2=0$, $\sum_{i=1}^{5} z_i \overline{z_i}=1$, $z_i \in {\Bbb C}$ ($1 \leq i \leq 5$),
which is the Brieskorn description of one of Milnor's exotic $7$-dimensional spheres. In fact, it is an
exotic $7$-sphere representing Milnor's standard generator of $\Theta_7$ (see \cite{Br1}, \cite{Br2}, \cite{Br3},
\cite{Br5}, \cite{Hi1}, \cite{Hi2}, \cite{Hi4}, \cite{HM}, \cite{Mi1}, \cite{Mi2}, \cite{Mi3} and \cite{Mi4}).

  We will study the following problems which are closely related to the $E_8$-singularity or more generally,
the exceptional singularities.

{\it Problem 1}. {\it The mysterious relation between exotic spheres, the icosahedron and $E_8$.}

  In his ICM 1970 talk \cite{Br3}, Brieskorn showed how to construct the singularity of type $ADE$ directly
from the simple complex Lie group of the same type. At the end of that paper \cite{Br3} Brieskorn
says:``Thus we see that there is a relation between exotic spheres, the icosahedron and $E_8$. But I still
do not understand why the regular polyhedra come in.'' In fact, even today there is some mystery in
these connections of such different parts of mathematics (see \cite{Gr}).

{\it Problem 2}. {\it Arnold's $A, D, E$ problem, especially for the appearance of the icosahedron.}

  In \cite{Ar2} and \cite{Ar6}, Arnold posed the following $A, D, E$ problem: ``The Dynkin diagrams $A_k$,
$D_k$, $E_6$, $E_7$, $E_8$ appear unexpectedly at the solution of so different classification problems, as
the classifications (1) of critical points of functions, (2) of regular polyhedra in ${\Bbb R}^3$, (3) the
category of linear spaces, (4) caustics, (5) wave fronts, (6) reflection groups, and (7) simple Lie groups. Some
connections between these objects are known, but in the majority of cases (e.g. in the cases $(1) \Leftrightarrow
(2) \Leftrightarrow (3)$) the coincidence of the answers to different problems has until now no explanation at all.
The problem ``$A, D, E$'' consists in finding such a general classification theorem from which the solutions of all listed
problems could be derived such that one could derive, let us say, the classification of simple stable singularities
of caustics from the classification of regular polyhedra.'' Moreover, in \cite{Ar3} and \cite{Ar5}, Arnold pointed
out that the link between the theory of singularities and the classification of regular polyhedra in three-dimensional
Euclidean space is rather mysteriously, especially for the appearance of the icosahedron.

{\it Problem 3}. {\it The relation between simple Lie algebras and Platonic solids.}

  At the end of \cite{Mc1}, McKay posed the following: ``Would not the Greeks appreciate the result that the simple Lie
algebras may be derived from the Platonic solids?''

{\it Problem 4}. {\it The Fermat-Catalan conjecture, especially the search for an eleventh solution.}

  The analytic construction of solutions of certain natural Diophantine equations is a problem of
central importance in number theory. In particular, it plays a central role in the study of the following
conjecture which contains within it essentially Fermat's Last Theorem and also the Catalan conjecture (see
\cite{DG} and \cite{D}).

{\bf Fermat-Catalan Conjecture}. {\it If $\frac{1}{p}+\frac{1}{q}+\frac{1}{r}<1$, then the Diophantine
equation $x^p+y^q=z^r$ in $x, y, z \in {\Bbb Z}$ with $\gcd(x, y, z)=1$, $xyz \neq 0$ and
$p, q, r$ are positive integers has no solutions except the following$:$
$$1+2^3=3^2, 2^5+7^2=3^4, 7^3+13^2=2^9, 2^7+17^3=71^2, 3^5+11^4=122^2,$$
$$17^7+76271^3=21063928^2, 1414^3+2213459^2=65^7, 9262^3+15312283^2=113^7,$$
$$43^8+96222^3=30042907^2, 33^8+1549034^2=15613^3.$$}

  There is a general approach to studying the equation $x^p+y^q=z^r$, which can be viewed as a kind of
nonabelian descent (see \cite{PSS}). The principal homogeneous spaces which arise in this descent are
algebraic curves equipped with a covering map to the projective line which is unramified outside of
$\{ 0, 1, \infty \}$ and is of signature $(p, q, r)$. This just means that the ramification indices of
the covering at the points above $0$, $1$ and $\infty$ divide $p$, $q$ and $r$ respectively. By a
beautiful application of ideas from classical invariant theory and modular forms, the analytic
construction of solutions were given in \cite{Be} and \cite{E} for $(p, q, r)=(2, 3, 3)$, $(2, 3, 4)$, $(2, 3, 5)$,
and in \cite{PSS} for $(p, q, r)=(2, 3, 7)$ which is the exceptional singularity $E_{12}$ (see \cite{Br4} and
\cite{BPR}). The integer solutions were then obtained from analytic construction of solutions (see \cite{Be},
\cite{E} and \cite{PSS}). A good starting point for the Fermat-Catalan Conjecture is to search for the eleventh
solution.

  In the present paper, we establish the invariant theory for $\text{PSL}(2, 13)$. Let us begin with the six-dimensional
representation of the finite simple group $\text{PSL}(2, 13)$ of order $1092$, which acts on the five-dimensional
projective space ${\Bbb P}^5=\{ (z_1, z_2, z_3, z_4, z_5, z_6): z_i \in {\Bbb C} \quad (i=1, 2, 3, 4, 5, 6) \}$.
This representation is defined over the cyclotomic field ${\Bbb Q}(e^{\frac{2 \pi i}{13}})$. Put
$$S=-\frac{1}{\sqrt{13}} \left(\matrix
  \zeta^{12}-\zeta & \zeta^{10}-\zeta^3 & \zeta^4-\zeta^9 & \zeta^5-\zeta^8 & \zeta^2-\zeta^{11} & \zeta^6-\zeta^7\\
  \zeta^{10}-\zeta^3 & \zeta^4-\zeta^9 & \zeta^{12}-\zeta & \zeta^2-\zeta^{11} & \zeta^6-\zeta^7 & \zeta^5-\zeta^8\\
  \zeta^4-\zeta^9 & \zeta^{12}-\zeta & \zeta^{10}-\zeta^3 & \zeta^6-\zeta^7 & \zeta^5-\zeta^8 & \zeta^2-\zeta^{11}\\
  \zeta^5-\zeta^8 & \zeta^2-\zeta^{11} & \zeta^6-\zeta^7 & \zeta-\zeta^{12} & \zeta^3-\zeta^{10} & \zeta^9-\zeta^4\\
  \zeta^2-\zeta^{11} & \zeta^6-\zeta^7 & \zeta^5-\zeta^8 & \zeta^3-\zeta^{10} & \zeta^9-\zeta^4 & \zeta-\zeta^{12}\\
  \zeta^6-\zeta^7 & \zeta^5-\zeta^8 & \zeta^2-\zeta^{11} & \zeta^9-\zeta^4 & \zeta-\zeta^{12} & \zeta^3-\zeta^{10}
\endmatrix\right)$$
and $T=\text{diag}(\zeta^7, \zeta^{11}, \zeta^8, \zeta^6, \zeta^2, \zeta^5)$ where $\zeta=\exp(2 \pi i/13)$.
We have $S^2=T^{13}=(ST)^3=1$. Let $G=\langle S, T \rangle$, then $G \cong \text{PSL}(2, 13)$. We construct some
$G$-invariant polynomials in six variables $z_1, \cdots, z_6$ which come from the Jacobian equation of degree
fourteen and the exotic equation of degree fourteen. Let
$$w_{\infty}=13 {\Bbb A}_0^2, \quad
  w_{\nu}=({\Bbb A}_0+\zeta^{\nu} {\Bbb A}_1+\zeta^{4 \nu} {\Bbb A}_2+\zeta^{9 \nu} {\Bbb A}_3+\zeta^{3 \nu}
  {\Bbb A}_4+\zeta^{12 \nu} {\Bbb A}_5+\zeta^{10 \nu} {\Bbb A}_6)^2\eqno{(1.2)}$$
for $\nu=0, 1, \cdots, 12$, where the senary quadratic forms (quadratic forms in six variables)
${\Bbb A}_j$ $(j=0, 1, \cdots, 6)$ are given by
$$\left\{\aligned
  {\Bbb A}_0 &=z_1 z_4+z_2 z_5+z_3 z_6,\\
  {\Bbb A}_1 &=z_1^2-2 z_3 z_4,\\
  {\Bbb A}_2 &=-z_5^2-2 z_2 z_4,\\
  {\Bbb A}_3 &=z_2^2-2 z_1 z_5,\\
  {\Bbb A}_4 &=z_3^2-2 z_2 z_6,\\
  {\Bbb A}_5 &=-z_4^2-2 z_1 z_6,\\
  {\Bbb A}_6 &=-z_6^2-2 z_3 z_5.
\endaligned\right.\eqno{(1.3)}$$
Then $w_{\infty}$, $w_{\nu}$ for $\nu=0, \cdots, 12$ form an algebraic equation of degree fourteen,
which is just the Jacobian equation of degree fourteen, whose roots are these $w_{\nu}$ and $w_{\infty}$.
On the other hand, set
$$\delta_{\infty}=13^2 {\Bbb G}_0, \quad
  \delta_{\nu}=-13 {\Bbb G}_0+\zeta^{\nu} {\Bbb G}_1+\zeta^{2 \nu} {\Bbb G}_2+\cdots+\zeta^{12 \nu} {\Bbb G}_{12}\eqno{(1.4)}$$
for $\nu=0, 1, \cdots, 12$, where the senary sextic forms (i.e., sextic forms in six variables) ${\Bbb G}_j$ $(j=0,
1, \cdots, 12)$ are given by
$$\left\{\aligned
  {\Bbb G}_0 &={\Bbb D}_0^2+{\Bbb D}_{\infty}^2,\\
  {\Bbb G}_1 &=-{\Bbb D}_7^2+2 {\Bbb D}_0 {\Bbb D}_1+10 {\Bbb D}_{\infty} {\Bbb D}_1
               +2 {\Bbb D}_2 {\Bbb D}_{12}-2 {\Bbb D}_3 {\Bbb D}_{11}-4 {\Bbb D}_4 {\Bbb D}_{10}
               -2 {\Bbb D}_9 {\Bbb D}_5,\\
  {\Bbb G}_2 &=-2 {\Bbb D}_1^2-4 {\Bbb D}_0 {\Bbb D}_2+6 {\Bbb D}_{\infty} {\Bbb D}_2
               -2 {\Bbb D}_4 {\Bbb D}_{11}+2 {\Bbb D}_5 {\Bbb D}_{10}-2 {\Bbb D}_6 {\Bbb D}_9
               -2 {\Bbb D}_7 {\Bbb D}_8,\\
  {\Bbb G}_3 &=-{\Bbb D}_8^2+2 {\Bbb D}_0 {\Bbb D}_3+10 {\Bbb D}_{\infty} {\Bbb D}_3
               +2 {\Bbb D}_6 {\Bbb D}_{10}-2 {\Bbb D}_9 {\Bbb D}_7-4 {\Bbb D}_{12} {\Bbb D}_4
               -2 {\Bbb D}_1 {\Bbb D}_2,\\
  {\Bbb G}_4 &=-{\Bbb D}_2^2+10 {\Bbb D}_0 {\Bbb D}_4-2 {\Bbb D}_{\infty} {\Bbb D}_4
               +2 {\Bbb D}_5 {\Bbb D}_{12}-2 {\Bbb D}_9 {\Bbb D}_8-4 {\Bbb D}_1 {\Bbb D}_3
               -2 {\Bbb D}_{10} {\Bbb D}_7,\\
  {\Bbb G}_5 &=-2 {\Bbb D}_9^2-4 {\Bbb D}_0 {\Bbb D}_5+6 {\Bbb D}_{\infty} {\Bbb D}_5
               -2 {\Bbb D}_{10} {\Bbb D}_8+2 {\Bbb D}_6 {\Bbb D}_{12}-2 {\Bbb D}_2 {\Bbb D}_3
               -2 {\Bbb D}_{11} {\Bbb D}_7,\\
  {\Bbb G}_6 &=-2 {\Bbb D}_3^2-4 {\Bbb D}_0 {\Bbb D}_6+6 {\Bbb D}_{\infty} {\Bbb D}_6
               -2 {\Bbb D}_{12} {\Bbb D}_7+2 {\Bbb D}_2 {\Bbb D}_4-2 {\Bbb D}_5 {\Bbb D}_1
               -2 {\Bbb D}_8 {\Bbb D}_{11},\\
  {\Bbb G}_7 &=-2 {\Bbb D}_{10}^2+6 {\Bbb D}_0 {\Bbb D}_7+4 {\Bbb D}_{\infty} {\Bbb D}_7
               -2 {\Bbb D}_1 {\Bbb D}_6-2 {\Bbb D}_2 {\Bbb D}_5-2 {\Bbb D}_8 {\Bbb D}_{12}
               -2 {\Bbb D}_9 {\Bbb D}_{11},\\
  {\Bbb G}_8 &=-2 {\Bbb D}_4^2+6 {\Bbb D}_0 {\Bbb D}_8+4 {\Bbb D}_{\infty} {\Bbb D}_8
               -2 {\Bbb D}_3 {\Bbb D}_5-2 {\Bbb D}_6 {\Bbb D}_2-2 {\Bbb D}_{11} {\Bbb D}_{10}
               -2 {\Bbb D}_1 {\Bbb D}_7,\\
  {\Bbb G}_9 &=-{\Bbb D}_{11}^2+2 {\Bbb D}_0 {\Bbb D}_9+10 {\Bbb D}_{\infty} {\Bbb D}_9
               +2 {\Bbb D}_5 {\Bbb D}_4-2 {\Bbb D}_1 {\Bbb D}_8-4 {\Bbb D}_{10} {\Bbb D}_{12}
               -2 {\Bbb D}_3 {\Bbb D}_6,\\
  {\Bbb G}_{10} &=-{\Bbb D}_5^2+10 {\Bbb D}_0 {\Bbb D}_{10}-2 {\Bbb D}_{\infty} {\Bbb D}_{10}
               +2 {\Bbb D}_6 {\Bbb D}_4-2 {\Bbb D}_3 {\Bbb D}_7-4 {\Bbb D}_9 {\Bbb D}_1
               -2 {\Bbb D}_{12} {\Bbb D}_{11},\\
  {\Bbb G}_{11} &=-2 {\Bbb D}_{12}^2+6 {\Bbb D}_0 {\Bbb D}_{11}+4 {\Bbb D}_{\infty} {\Bbb D}_{11}
               -2 {\Bbb D}_9 {\Bbb D}_2-2 {\Bbb D}_5 {\Bbb D}_6-2 {\Bbb D}_7 {\Bbb D}_4
               -2 {\Bbb D}_3 {\Bbb D}_8,\\
  {\Bbb G}_{12} &=-{\Bbb D}_6^2+10 {\Bbb D}_0 {\Bbb D}_{12}-2 {\Bbb D}_{\infty} {\Bbb D}_{12}
               +2 {\Bbb D}_2 {\Bbb D}_{10}-2 {\Bbb D}_1 {\Bbb D}_{11}-4 {\Bbb D}_3 {\Bbb D}_9
               -2 {\Bbb D}_4 {\Bbb D}_8.
\endaligned\right.\eqno{(1.5)}$$
Here, the senary cubic forms (cubic forms in six variables) ${\Bbb D}_j$ $(j=0, 1, \cdots, 12, \infty)$
are given as follows:
$$\left\{\aligned
  {\Bbb D}_0 &=z_1 z_2 z_3,\\
  {\Bbb D}_1 &=2 z_2 z_3^2+z_2^2 z_6-z_4^2 z_5+z_1 z_5 z_6,\\
  {\Bbb D}_2 &=-z_6^3+z_2^2 z_4-2 z_2 z_5^2+z_1 z_4 z_5+3 z_3 z_5 z_6,\\
  {\Bbb D}_3 &=2 z_1 z_2^2+z_1^2 z_5-z_4 z_6^2+z_3 z_4 z_5,\\
  {\Bbb D}_4 &=-z_2^2 z_3+z_1 z_6^2-2 z_4^2 z_6-z_1 z_3 z_5,\\
  {\Bbb D}_5 &=-z_4^3+z_3^2 z_5-2 z_3 z_6^2+z_2 z_5 z_6+3 z_1 z_4 z_6,\\
  {\Bbb D}_6 &=-z_5^3+z_1^2 z_6-2 z_1 z_4^2+z_3 z_4 z_6+3 z_2 z_4 z_5,\\
  {\Bbb D}_7 &=-z_2^3+z_3 z_4^2-z_1 z_3 z_6-3 z_1 z_2 z_5+2 z_1^2 z_4,\\
  {\Bbb D}_8 &=-z_1^3+z_2 z_6^2-z_2 z_3 z_5-3 z_1 z_3 z_4+2 z_3^2 z_6,\\
  {\Bbb D}_9 &=2 z_1^2 z_3+z_3^2 z_4-z_5^2 z_6+z_2 z_4 z_6,\\
  {\Bbb D}_{10} &=-z_1 z_3^2+z_2 z_4^2-2 z_4 z_5^2-z_1 z_2 z_6,\\
  {\Bbb D}_{11} &=-z_3^3+z_1 z_5^2-z_1 z_2 z_4-3 z_2 z_3 z_6+2 z_2^2 z_5,\\
  {\Bbb D}_{12} &=-z_1^2 z_2+z_3 z_5^2-2 z_5 z_6^2-z_2 z_3 z_4,\\
  {\Bbb D}_{\infty}&=z_4 z_5 z_6.
\endaligned\right.\eqno{(1.6)}$$
Then $\delta_{\infty}$, $\delta_{\nu}$ for $\nu=0, \cdots, 12$ form an algebraic equation of degree fourteen
which is not the Jacobian equation of degree fourteen. We call it the exotic equation of degree fourteen. Let
$$\Phi_{12}=-\frac{1}{13 \cdot 52} \left(\sum_{\nu=0}^{12} \delta_{\nu}^2+\delta_{\infty}^2\right), \quad
  \Phi_{18}=\frac{1}{13 \cdot 6} \left(\sum_{\nu=0}^{12} \delta_{\nu}^3+\delta_{\infty}^3\right),\eqno{(1.7)}$$
$$\Phi_{20}=\frac{1}{13 \cdot 25} \left(\sum_{\nu=0}^{12} w_{\nu}^5+w_{\infty}^5\right), \quad
  \Phi_{30}=-\frac{1}{13 \cdot 1315} \left(\sum_{\nu=0}^{12} \delta_{\nu}^5+\delta_{\infty}^5\right),\eqno{(1.8)}$$
and $x_i(z)=\eta(z) a_i(z)$ $(1 \leq i \leq 6)$, where
$$\left\{\aligned
  a_1(z) &:=e^{-\frac{11 \pi i}{26}} \theta \left[\matrix \frac{11}{13}\\ 1 \endmatrix\right](0, 13z),\\
  a_2(z) &:=e^{-\frac{7 \pi i}{26}} \theta \left[\matrix \frac{7}{13}\\ 1 \endmatrix\right](0, 13z),\\
  a_3(z) &:=e^{-\frac{5 \pi i}{26}} \theta \left[\matrix \frac{5}{13}\\ 1 \endmatrix\right](0, 13z),\\
  a_4(z) &:=-e^{-\frac{3 \pi i}{26}} \theta \left[\matrix \frac{3}{13}\\ 1 \endmatrix\right](0, 13z),\\
  a_5(z) &:=e^{-\frac{9 \pi i}{26}} \theta \left[\matrix \frac{9}{13}\\ 1 \endmatrix\right](0, 13z),\\
  a_6(z) &:=e^{-\frac{\pi i}{26}} \theta \left[\matrix \frac{1}{13}\\ 1 \endmatrix\right](0, 13z)
\endaligned\right.\eqno{(1.9)}$$
are theta constants of order $13$ and $\eta(z):=q^{\frac{1}{24}} \prod_{n=1}^{\infty} (1-q^n)$ with
$q=e^{2 \pi i z}$ is the Dedekind eta function which are all defined in the upper-half plane
${\Bbb H}=\{ z \in {\Bbb C}: \text{Im}(z)>0 \}$. Our main result is to show that there exists at least
two kinds of symmetry groups associated with the equation of $E_8$-singularity: one is the icosahedral
group (the standard structure), the other is the group $\text{PSL}(2, 13)$ (the exotic structure).

{\bf Theorem 1.1 (Main Theorem 1)}. {\it The exotic structure on the $E_8$-singularity is given by
the following equations:
$$\Phi_{20}^3-\Phi_{30}^2=1728 \Phi_{12}^5, \quad \Phi_{20}^3-\Phi_{12}^2 \Phi_{18}^2=1728 \Phi_{12}^5,\eqno{(1.10)}$$
where we use the abbreviation $\Phi_j=\Phi_j(x_1(z), \cdots, x_6(z))$ for $j=12, 18, 20$ and $30$.
As polynomials in six variables $z_1, \cdots, z_6$, $\Phi_{12}$, $\Phi_{18}$, $\Phi_{20}$ and
$\Phi_{30}$ are $G$-invariant polynomials.}

  As a consequence, Theorem 1.1 shows that $E_8$ is not uniquely determined by the icosahedron. This gives
a negative answer to Problem 1, i.e., the icosahedron does not necessarily appear in the triple (exotic
spheres, icosahedron, $E_8$). The group $\text{PSL}(2, 13)$ can take its place and there is the other
triple (exotic spheres, $\text{PSL}(2, 13)$, $E_8$). The higher dimensional liftings of these two distinct
symmetry groups and modular interpretations on the equation of $E_8$-singularity give the same Milnor's standard
generator of $\Theta_7$. Simultaneously, Theorem 1.1 gives a counterexample to Problem 2. It shows that the equivalence
between critical points of functions and regular polyhedra in ${\Bbb R}^3$ in Problem 2 is not true. In particular,
the appearance of the icosahedron is not necessary. The group $\text{PSL}(2, 13)$ can take its place, which is the
symmetry group of the regular maps $\{ 7, 3 \}_{13}$, $\{ 13, 3 \}_7$ or $\{ 13, 7 \}_3$ (see \cite{CoM}, p.140).
Theorem 1.1 also shows that the simple Lie algebra $E_8$ may not be necessarily derived from the icosahedron,
which gives a negative answer to Problem 3. On the other hand, Theorem 1.1 provides an analytic construction of
solutions for the Diophantine equation $x^p+y^q=z^r$ in the case $(p, q, r)=(2, 3, 5)$ which is different from
the solutions given in \cite{E}.

  Let
$$\Phi_{32}=-\frac{1}{13 \cdot 1840} \left(\sum_{\nu=0}^{12} w_{\nu}^8+w_{\infty}^8\right), \quad
  \Phi_{42}=\frac{1}{13 \cdot 226842} \left(\sum_{\nu=0}^{12} \delta_{\nu}^7+\delta_{\infty}^7\right),\eqno{(1.11)}$$
$$\Phi_{44}=\frac{1}{13 \cdot 146905} \left(\sum_{\nu=0}^{12} w_{\nu}^{11}+w_{\infty}^{11}\right),\eqno{(1.12)}$$

{\bf Theorem 1.2. (Main Theorem 2)}. {\it The exceptional singularities
$Q_{18}: x^3+y^8+y z^2=0$, $E_{20}: x^3+y^{11}+z^2=0$ and $x^7+x^2 y^3+z^2=0$
can be endowed with the symmetry group $G$ and have the following modular parametrizations:
$$\left\{\aligned
  \Phi_{20}^3 \Phi_{12}^2-\Phi_{42}^2 &=1728 \Phi_{12}^7,\\
  \Phi_{32}^3-\Phi_{12}^5 \Phi_{18}^2 &=1728 \Phi_{12}^8,\\
  \Phi_{32}^3-\Phi_{12}^3 \Phi_{30}^2 &=1728 \Phi_{12}^8,\\
    \Phi_{32}^3-\Phi_{12} \Phi_{42}^2 &=1728 \Phi_{12}^8,\\
  \Phi_{44}^3-\Phi_{12}^8 \Phi_{18}^2 &=1728 \Phi_{12}^{11},\\
  \Phi_{44}^3-\Phi_{12}^6 \Phi_{30}^2 &=1728 \Phi_{12}^{11},\\
  \Phi_{44}^3-\Phi_{12}^4 \Phi_{42}^2 &=1728 \Phi_{12}^{11},
\endaligned\right.\eqno{(1.13)}$$
where we use the abbreviation $\Phi_j=\Phi_j(x_1(z), \cdots, x_6(z))$ for $j=32, 42$ and $44$.
As polynomials in six variables $z_1, \cdots, z_6$, $\Phi_{12}$, $\Phi_{18}$, $\Phi_{20}$, $\Phi_{30}$
$\Phi_{32}$, $\Phi_{42}$ and $\Phi_{44}$ are $G$-invariant polynomials.}

{\bf Corollary 1.3}. {\it The Diophantine equation $x^2+y^3=z^{11}$ has the following analytic
construction of solutions:
$$\left\{\aligned
  \Phi_{44}^3-\Phi_{12}^8 \Phi_{18}^2 &=1728 \Phi_{12}^{11},\\
  \Phi_{44}^3-\Phi_{12}^6 \Phi_{30}^2 &=1728 \Phi_{12}^{11},\\
  \Phi_{44}^3-\Phi_{12}^4 \Phi_{42}^2 &=1728 \Phi_{12}^{11},
\endaligned\right.\eqno{(1.14)}$$
where we use the abbreviation $\Phi_j=\Phi_j(x_1(z), \cdots, x_6(z))$ for $j=12, 18, 30, 42$ and $44$.}

  We can multiply the Diophantine equation $x^2+y^3=z^{11}$ with $\gcd(x, y, z)=1$ by $3^{36} 2^6$ to obtain
$(-3^{12} 2^2 y)^3-(3^{18} 2^3 x)^2=1728 (-3^3 z)^{11}$. Hence, Corollary 1.3 provides a new analytic construction
of solutions for the Diophantine equation $x^p+y^q=z^r$ in the more difficult case $(p, q, r)=(2, 3, 11)$, which
goes beyond the known ten solutions in the Fermat-Catalan Conjecture and gives an answer to Problem 4. Theorem 1.1
and Theorem 1.2 give an answer to a problem dating back to the works of Klein (see the end of section four).

  In number theory and arithmetical algebraic geometry, a central problem is the modularity of algebraic
varieties. Our example of exotic structure on the equation of $E_8$-singularity shows
that the algebraic structure of a variety does not completely determine its arithmetical (modular) structure.
This leads to the study of the relationships between algebraic property and arithmetical (modular) property
on an algebraic variety. Note that the Klein quartic curve has two kinds of modular parametrizations with the
same level and the same symmetry group. One is the modular curve $X(7)$ of level $7$ with symmetry group
$\Gamma/\Gamma(7) \cong \text{PSL}(2, 7)$. The other is the Shimura curve ${\Bbb H}/\Gamma({\frak o}, {\frak p})$
for a certain quaternion algebra ${\frak o}$ over a cubic field ${\Bbb Q}(\cos \frac{2 \pi}{7})$ and
$N({\frak p})=7$, whose symmetry group is $\Gamma({\frak o})/\Gamma({\frak o}, {\frak p}) \cong \text{PSL}(2, 7)$
(see \cite{Sh}). In contrast with it, our example has different levels and different symmetry groups.

  In fact, $G$ is a finite subgroup of $\text{SL}(6, {\Bbb C})$ (see \cite{L}) and ${\Bbb C}^6/G$ is a
six-dimensional quotient singularity. Theorem 1.1 and Theorem 1.2 imply that each of the following nine
triples $(\Phi_{12}, \Phi_{18}, \Phi_{20})$, $(\Phi_{12}, \Phi_{20}, \Phi_{30})$,
$(\Phi_{12}, \Phi_{20}, \Phi_{42})$, $(\Phi_{12}, \Phi_{18}, \Phi_{32})$, $(\Phi_{12}, \Phi_{30}, \Phi_{32})$,
$(\Phi_{12}, \Phi_{32}, \Phi_{42})$, $(\Phi_{12}, \Phi_{18}, \Phi_{44})$, $(\Phi_{12}, \Phi_{30}, \Phi_{44})$
and $(\Phi_{12}, \Phi_{42}, \Phi_{44})$ gives a map from ${\Bbb C}^6/G$ to ${\Bbb C}^3$ with the corresponding
syzygy which relates the three polynomials in six variables when evaluating at theta constants of order thirteen.

  This paper consists of four sections. In section two, we give the standard structure on the
$E_8$-singularity. In section three, we obtain the invariant theory for $\text{PSL}(2, 13)$.
In particular, we construct the senary quadratic forms ${\Bbb A}_j$ ($0 \leq j \leq 6$),  the
senary cubic forms ${\Bbb D}_j$ ($j=0, 1, \cdots, 12, \infty$) and the senary sextic forms
${\Bbb G}_j$ ($0 \leq j \leq 12$). From ${\Bbb A}_j$ we construct the Jacobian equation of degree
fourteen. From ${\Bbb D}_j$ and ${\Bbb G}_j$ we construct the exotic equation of degree fourteen.
It should be pointed out that the appearance of exotic equation is a new phenomenon which does not
appear in the invariant theory for $\text{PSL}(2, 5)$, $\text{PSL}(2, 7)$ and $\text{PSL}(2, 11)$
(see \cite{K}, \cite{K1}, \cite{K2}, \cite{K3}, \cite{K4}, \cite{KF1} and \cite{KF2}). Combining Jacobian equation
with exotic equation, we obtain the polynomials $\Phi_{12}$, $\Phi_{18}$, $\Phi_{20}$ and $\Phi_{30}$
which are invariant under the action of $\text{PSL}(2, 13)$. Together with theta constants of order
thirteen, we give the exotic structure on the $E_8$-singularity and prove Theorem 1.1. In section
four, we obtain the invariant polynomials $\Phi_{32}$, $\Phi_{42}$, and $\Phi_{44}$. Together with
theta constants of order thirteen, we give modular parametrizations for the exceptional singularities
$Q_{18}$, $E_{20}$ and $x^7+x^2 y^3+z^2=0$ and prove Theorem 1.2.

\vskip 0.5 cm

\centerline{\bf 2. Standard structure on the $E_8$-singularity: the icosahedron}

\vskip 0.5 cm

  In this section, we will give the standard structure on the equation of $E_8$-singularity: the
classical icosahedral equation with icosahedral symmetry, which has a modular interpretation in
terms of theta constants of order five.

  The most famous source on icosahedral symmetry is the celebrated book \cite{K} of Felix Klein
on the icosahedron and the solution of quintic equations which appeared for the first time in
1884. Its main objective is to show that the solution of general quintic equations can be
reduced to that of locally inverting certain quotients of actions of the icosahedral group
$G \cong A_5$ on spaces related to the usual geometric realization of the regular icosahedron in
three-dimensional space. In particular, Klein (see \cite{K}) gave a parametric solution
of the $E_8$-singularity using the icosahedron by homogeneous polynomials $T$, $H$, $f$ in two
variables of degrees $30$, $20$, $12$ with integral coefficients. Here,
$$f=z_1 z_2 (z_1^{10}+11 z_1^5 z_2^5-z_2^{10}),$$
$$H=\frac{1}{121} \vmatrix \format \c \quad & \c\\
    \frac{\partial^2 f}{\partial z_1^2} &
    \frac{\partial^2 f}{\partial z_1 \partial z_2}\\
    \frac{\partial^2 f}{\partial z_2 \partial z_1} &
    \frac{\partial^2 f}{\partial z_2^2}
    \endvmatrix
  =-(z_1^{20}+z_2^{20})+228 (z_1^{15} z_2^5-z_1^5 z_2^{15})-494 z_1^{10} z_2^{10},$$
$$T=-\frac{1}{20} \vmatrix \format \c \quad & \c\\
    \frac{\partial f}{\partial z_1} &
    \frac{\partial f}{\partial z_2}\\
    \frac{\partial H}{\partial z_1} &
    \frac{\partial H}{\partial z_2}
    \endvmatrix
  =(z_1^{30}+z_2^{30})+522 (z_1^{25} z_2^5-z_1^5 z_2^{25})-10005 (z_1^{20} z_2^{10}+z_1^{10} z_2^{20}).$$
They satisfy the famous icosahedral equation
$$T^2+H^3=1728 f^5.\eqno{(2.1)}$$
In fact, $f$, $H$ and $T$ are invariant polynomials under the action of the icosahedral group.
Essentially the same relation had been found a few years earlier by Schwarz (see \cite{Sch}), who
considered three polynomials $\varphi_{12}$, $\varphi_{20}$ and $\varphi_{30}$ whose roots correspond
to the vertices, the midpoints of the faces and the midpoints of the edges of an icosahedron inscribed
in the Riemann sphere. He obtained the identity $\varphi_{20}^3-1728 \varphi_{12}^5=\varphi_{30}^2$.
Thus we see that from the very beginning there was a close relation between the $E_8$-singularity and the
icosahedron. Moreover, the icosahedral equation (2.1) can be interpreted in terms of modular forms which
was also known by Klein (see \cite{KF1}, p. 631). Let $x_1(z)=\eta(z) a(z)$ and $x_2(z)=\eta(z) b(z)$, where
$$a(z)=e^{-\frac{3 \pi i}{10}} \theta \left[\matrix
       \frac{3}{5}\\ 1 \endmatrix\right](0, 5z), \quad
  b(z)=e^{-\frac{\pi i}{10}} \theta \left[\matrix
       \frac{1}{5}\\ 1 \endmatrix\right](0, 5z)$$
are theta constants of order five and $\eta(z):=q^{\frac{1}{24}} \prod_{n=1}^{\infty} (1-q^n)$ with
$q=e^{2 \pi i z}$ is the Dedekind eta function which are all defined in the upper-half plane
${\Bbb H}=\{ z \in {\Bbb C}: \text{Im}(z)>0 \}$. Then
$$\left\{\aligned
  f(x_1(z), x_2(z)) &=-\Delta(z),\\
  H(x_1(z), x_2(z)) &=-\eta(z)^8 \Delta(z)E_4(z),\\
  T(x_1(z), x_2(z)) &=\Delta(z)^2 E_6(z),
\endaligned\right.$$
where $E_4(z):=\frac{1}{2} \sum_{m, n \in {\Bbb Z}, (m, n)=1} \frac{1}{(mz+n)^4}$ and
$E_6(z):=\frac{1}{2} \sum_{m, n \in {\Bbb Z}, (m, n)=1} \frac{1}{(mz+n)^6}$ are Eisenstein series
of weight $4$ and $6$, and $\Delta(z)=\eta(z)^{24}$ is the discriminant. The relations
$$j(z):=\frac{E_4(z)^3}{\Delta(z)}=\frac{H(x_1(z), x_2(z))^3}{f(x_1(z), x_2(z))^5}, \quad
  j(z)-1728=\frac{E_6(z)^2}{\Delta(z)}=-\frac{T(x_1(z), x_2(z))^2}{f(x_1(z), x_2(z))^5}$$
give the icosahedral equation (2.1) in terms of theta constants of order five.

\vskip 0.5 cm

\centerline{\bf 3. Exotic structure on the $E_8$-singularity: $\text{PSL}(2, 13)$}

\vskip 0.5 cm

  In this section, we will give the exotic structure on the equation of $E_8$-singularity: the
symmetry group is the simple group $\text{PSL}(2, 13)$ and the equation has a modular interpretation
in terms of theta constants of order thirteen.

  At first, we will study the six-dimensional representation of the finite simple group $\text{PSL}(2, 13)$
of order $1092$, which acts on the five-dimensional projective space
${\Bbb P}^5=\{ (z_1, z_2, z_3, z_4, z_5, z_6): z_i \in {\Bbb C} \quad (i=1, 2, 3, 4, 5, 6) \}$.
This representation is defined over the cyclotomic field ${\Bbb Q}(e^{\frac{2 \pi i}{13}})$. Put
$$S=-\frac{1}{\sqrt{13}} \left(\matrix
  \zeta^{12}-\zeta & \zeta^{10}-\zeta^3 & \zeta^4-\zeta^9 & \zeta^5-\zeta^8 & \zeta^2-\zeta^{11} & \zeta^6-\zeta^7\\
  \zeta^{10}-\zeta^3 & \zeta^4-\zeta^9 & \zeta^{12}-\zeta & \zeta^2-\zeta^{11} & \zeta^6-\zeta^7 & \zeta^5-\zeta^8\\
  \zeta^4-\zeta^9 & \zeta^{12}-\zeta & \zeta^{10}-\zeta^3 & \zeta^6-\zeta^7 & \zeta^5-\zeta^8 & \zeta^2-\zeta^{11}\\
  \zeta^5-\zeta^8 & \zeta^2-\zeta^{11} & \zeta^6-\zeta^7 & \zeta-\zeta^{12} & \zeta^3-\zeta^{10} & \zeta^9-\zeta^4\\
  \zeta^2-\zeta^{11} & \zeta^6-\zeta^7 & \zeta^5-\zeta^8 & \zeta^3-\zeta^{10} & \zeta^9-\zeta^4 & \zeta-\zeta^{12}\\
  \zeta^6-\zeta^7 & \zeta^5-\zeta^8 & \zeta^2-\zeta^{11} & \zeta^9-\zeta^4 & \zeta-\zeta^{12} & \zeta^3-\zeta^{10}
\endmatrix\right)\eqno{(3.1)}$$
and
$$T=\text{diag}(\zeta^7, \zeta^{11}, \zeta^8, \zeta^6, \zeta^2, \zeta^5),\eqno{(3.2)}$$
where $\zeta=\exp(2 \pi i/13)$. We have
$$S^2=T^{13}=(ST)^3=1.\eqno{(3.3)}$$
Let $G=\langle S, T \rangle$, then $G \cong \text{PSL}(2, 13)$ (see \cite{Y1}, Theorem 3.1).

  Put $\theta_1=\zeta+\zeta^3+\zeta^9$, $\theta_2=\zeta^2+\zeta^6+\zeta^5$, $\theta_3=\zeta^4+\zeta^{12}+\zeta^{10}$,
and $\theta_4=\zeta^8+\zeta^{11}+\zeta^7$. We find that
$$\left\{\aligned
  &\theta_1+\theta_2+\theta_3+\theta_4=-1,\\
  &\theta_1 \theta_2+\theta_1 \theta_3+\theta_1 \theta_4+\theta_2 \theta_3+\theta_2 \theta_4+\theta_3 \theta_4=2,\\
  &\theta_1 \theta_2 \theta_3+\theta_1 \theta_2 \theta_4+\theta_1 \theta_3 \theta_4+\theta_2 \theta_3 \theta_4=4,\\
  &\theta_1 \theta_2 \theta_3 \theta_4=3.
\endaligned\right.$$
Hence, $\theta_1$, $\theta_2$, $\theta_3$ and $\theta_4$ satisfy
the quartic equation $z^4+z^3+2 z^2-4z+3=0$,
which can be decomposed as two quadratic equations
$$\left(z^2+\frac{1+\sqrt{13}}{2} z+\frac{5+\sqrt{13}}{2}\right)
  \left(z^2+\frac{1-\sqrt{13}}{2} z+\frac{5-\sqrt{13}}{2}\right)=0$$
over the real quadratic field ${\Bbb Q}(\sqrt{13})$. Therefore, the
four roots are given as follows:
$$\left\{\aligned
  \theta_1=\frac{1}{4} \left(-1+\sqrt{13}+\sqrt{-26+6 \sqrt{13}}\right),\\
  \theta_2=\frac{1}{4} \left(-1-\sqrt{13}+\sqrt{-26-6 \sqrt{13}}\right),\\
  \theta_3=\frac{1}{4} \left(-1+\sqrt{13}-\sqrt{-26+6 \sqrt{13}}\right),\\
  \theta_4=\frac{1}{4} \left(-1-\sqrt{13}-\sqrt{-26-6 \sqrt{13}}\right).
\endaligned\right.$$
Moreover, we find that
$$\left\{\aligned
  \theta_1+\theta_3+\theta_2+\theta_4 &=-1,\\
  \theta_1+\theta_3-\theta_2-\theta_4 &=\sqrt{13},\\
  \theta_1-\theta_3-\theta_2+\theta_4 &=-\sqrt{-13+2 \sqrt{13}},\\
  \theta_1-\theta_3+\theta_2-\theta_4 &=\sqrt{-13-2 \sqrt{13}}.
\endaligned\right.$$

  Let us study the action of $S T^{\nu}$ on ${\Bbb P}^5$, where $\nu=0, 1, \cdots, 12$. Put
$$\alpha=\zeta+\zeta^{12}-\zeta^5-\zeta^8, \quad
   \beta=\zeta^3+\zeta^{10}-\zeta^2-\zeta^{11}, \quad
   \gamma=\zeta^9+\zeta^4-\zeta^6-\zeta^7.$$
We find that
$$\aligned
  &13 ST^{\nu}(z_1) \cdot ST^{\nu}(z_4)\\
=&\beta z_1 z_4+\gamma z_2 z_5+\alpha z_3 z_6+\\
 &+\gamma \zeta^{\nu} z_1^2+\alpha \zeta^{9 \nu} z_2^2+\beta \zeta^{3 \nu} z_3^2
  -\gamma \zeta^{12 \nu} z_4^2-\alpha \zeta^{4 \nu} z_5^2-\beta \zeta^{10 \nu} z_6^2+\\
 &+(\alpha-\beta) \zeta^{5 \nu} z_1 z_2+(\beta-\gamma) \zeta^{6 \nu} z_2 z_3
  +(\gamma-\alpha) \zeta^{2 \nu} z_1 z_3+\\
 &+(\beta-\alpha) \zeta^{8 \nu} z_4 z_5+(\gamma-\beta) \zeta^{7 \nu} z_5 z_6
  +(\alpha-\gamma) \zeta^{11 \nu} z_4 z_6+\\
 &-(\alpha+\beta) \zeta^{\nu} z_3 z_4-(\beta+\gamma) \zeta^{9 \nu} z_1 z_5
  -(\gamma+\alpha) \zeta^{3 \nu} z_2 z_6+\\
 &-(\alpha+\beta) \zeta^{12 \nu} z_1 z_6-(\beta+\gamma) \zeta^{4 \nu} z_2 z_4
  -(\gamma+\alpha) \zeta^{10 \nu} z_3 z_5.
\endaligned$$
$$\aligned
  &13 ST^{\nu}(z_2) \cdot ST^{\nu}(z_5)\\
=&\gamma z_1 z_4+\alpha z_2 z_5+\beta z_3 z_6+\\
 &+\alpha \zeta^{\nu} z_1^2+\beta \zeta^{9 \nu} z_2^2+\gamma \zeta^{3 \nu} z_3^2
  -\alpha \zeta^{12 \nu} z_4^2-\beta \zeta^{4 \nu} z_5^2-\gamma \zeta^{10 \nu} z_6^2+\\
 &+(\beta-\gamma) \zeta^{5 \nu} z_1 z_2+(\gamma-\alpha) \zeta^{6 \nu} z_2 z_3
  +(\alpha-\beta) \zeta^{2 \nu} z_1 z_3+\\
 &+(\gamma-\beta) \zeta^{8 \nu} z_4 z_5+(\alpha-\gamma) \zeta^{7 \nu} z_5 z_6
  +(\beta-\alpha) \zeta^{11 \nu} z_4 z_6+\\
 &-(\beta+\gamma) \zeta^{\nu} z_3 z_4-(\gamma+\alpha) \zeta^{9 \nu} z_1 z_5
  -(\alpha+\beta) \zeta^{3 \nu} z_2 z_6+\\
 &-(\beta+\gamma) \zeta^{12 \nu} z_1 z_6-(\gamma+\alpha) \zeta^{4 \nu} z_2 z_4
  -(\alpha+\beta) \zeta^{10 \nu} z_3 z_5.
\endaligned$$
$$\aligned
  &13 ST^{\nu}(z_3) \cdot ST^{\nu}(z_6)\\
=&\alpha z_1 z_4+\beta z_2 z_5+\gamma z_3 z_6+\\
 &+\beta \zeta^{\nu} z_1^2+\gamma \zeta^{9 \nu} z_2^2+\alpha \zeta^{3 \nu} z_3^2
  -\beta \zeta^{12 \nu} z_4^2-\gamma \zeta^{4 \nu} z_5^2-\alpha \zeta^{10 \nu} z_6^2+\\
 &+(\gamma-\alpha) \zeta^{5 \nu} z_1 z_2+(\alpha-\beta) \zeta^{6 \nu} z_2 z_3
  +(\beta-\gamma) \zeta^{2 \nu} z_1 z_3+\\
 &+(\alpha-\gamma) \zeta^{8 \nu} z_4 z_5+(\beta-\alpha) \zeta^{7 \nu} z_5 z_6
  +(\gamma-\beta) \zeta^{11 \nu} z_4 z_6+\\
 &-(\gamma+\alpha) \zeta^{\nu} z_3 z_4-(\alpha+\beta) \zeta^{9 \nu} z_1 z_5
  -(\beta+\gamma) \zeta^{3 \nu} z_2 z_6+\\
 &-(\gamma+\alpha) \zeta^{12 \nu} z_1 z_6-(\alpha+\beta) \zeta^{4 \nu} z_2 z_4
  -(\beta+\gamma) \zeta^{10 \nu} z_3 z_5.
\endaligned$$
Note that $\alpha+\beta+\gamma=\sqrt{13}$, we find that
$$\aligned
  &\sqrt{13} \left[ST^{\nu}(z_1) \cdot ST^{\nu}(z_4)+ST^{\nu}(z_2) \cdot ST^{\nu}(z_5)+ST^{\nu}(z_3) \cdot ST^{\nu}(z_6)\right]\\
 =&(z_1 z_4+z_2 z_5+z_3 z_6)+(\zeta^{\nu} z_1^2+\zeta^{9 \nu} z_2^2+\zeta^{3 \nu} z_3^2)
  -(\zeta^{12 \nu} z_4^2+\zeta^{4 \nu} z_5^2+\zeta^{10 \nu} z_6^2)+\\
  &-2(\zeta^{\nu} z_3 z_4+\zeta^{9 \nu} z_1 z_5+\zeta^{3 \nu} z_2 z_6)
   -2(\zeta^{12 \nu} z_1 z_6+\zeta^{4 \nu} z_2 z_4+\zeta^{10 \nu} z_3 z_5).
\endaligned$$
Let
$$\varphi_{\infty}(z_1, z_2, z_3, z_4, z_5, z_6)=\sqrt{13} (z_1 z_4+z_2 z_5+z_3 z_6)\eqno{(3.4)}$$
and
$$\varphi_{\nu}(z_1, z_2, z_3, z_4, z_5, z_6)=\varphi_{\infty}(ST^{\nu}(z_1, z_2, z_3, z_4, z_5, z_6))\eqno{(3.5)}$$
for $\nu=0, 1, \cdots, 12$. Then
$$\aligned
  \varphi_{\nu}
=&(z_1 z_4+z_2 z_5+z_3 z_6)+\zeta^{\nu} (z_1^2-2 z_3 z_4)+\zeta^{4 \nu} (-z_5^2-2 z_2 z_4)+\\
 &+\zeta^{9 \nu} (z_2^2-2 z_1 z_5)+\zeta^{3 \nu} (z_3^2-2 z_2 z_6)+
   \zeta^{12 \nu} (-z_4^2-2 z_1 z_6)+\zeta^{10 \nu} (-z_6^2-2 z_3 z_5).
\endaligned\eqno{(3.6)}$$
This leads us to define the following senary quadratic forms
(quadratic forms in six variables):
$$\left\{\aligned
  {\Bbb A}_0 &=z_1 z_4+z_2 z_5+z_3 z_6,\\
  {\Bbb A}_1 &=z_1^2-2 z_3 z_4,\\
  {\Bbb A}_2 &=-z_5^2-2 z_2 z_4,\\
  {\Bbb A}_3 &=z_2^2-2 z_1 z_5,\\
  {\Bbb A}_4 &=z_3^2-2 z_2 z_6,\\
  {\Bbb A}_5 &=-z_4^2-2 z_1 z_6,\\
  {\Bbb A}_6 &=-z_6^2-2 z_3 z_5.
\endaligned\right.\eqno{(3.7)}$$
Hence,
$$\sqrt{13} ST^{\nu}({\Bbb A}_0)={\Bbb A}_0+\zeta^{\nu} {\Bbb A}_1+\zeta^{4 \nu} {\Bbb A}_2+
  \zeta^{9 \nu} {\Bbb A}_3+\zeta^{3 \nu} {\Bbb A}_4+\zeta^{12 \nu} {\Bbb A}_5+\zeta^{10 \nu} {\Bbb A}_6.
  \eqno{(3.8)}$$
Let $H:=Q^5 P^2 \cdot P^2 Q^6 P^8 \cdot Q^5 P^2 \cdot P^3 Q$ where $P=S T^{-1} S$ and $Q=S T^3$.
Then (see \cite{Y2}, p.27)
$$H=\left(\matrix
  0 &  0 &  0 & 0 & 0 & 1\\
  0 &  0 &  0 & 1 & 0 & 0\\
  0 &  0 &  0 & 0 & 1 & 0\\
  0 &  0 & -1 & 0 & 0 & 0\\
 -1 &  0 &  0 & 0 & 0 & 0\\
  0 & -1 &  0 & 0 & 0 & 0
\endmatrix\right).\eqno{(3.9)}$$
Note that $H^6=1$ and $H^{-1} T H=-T^4$. Thus,
$\langle H, T \rangle \cong {\Bbb Z}_{13} \rtimes {\Bbb Z}_6$. Hence, it is a maximal subgroup of order $78$ of $G$ with
index $14$ (see \cite{CC}). We find that $\varphi_{\infty}^2$ is
invariant under the action of the maximal subgroup $\langle H, T
\rangle$. Note that
$$\varphi_{\infty}=\sqrt{13} {\Bbb A}_0, \quad
  \varphi_{\nu}={\Bbb A}_0+\zeta^{\nu} {\Bbb A}_1+\zeta^{4 \nu} {\Bbb A}_2+
  \zeta^{9 \nu} {\Bbb A}_3+\zeta^{3 \nu} {\Bbb A}_4+\zeta^{12 \nu} {\Bbb A}_5+\zeta^{10 \nu} {\Bbb A}_6$$
for $\nu=0, 1, \cdots, 12$. Let $w=\varphi^2$,
$w_{\infty}=\varphi_{\infty}^2$ and $w_{\nu}=\varphi_{\nu}^2$.
Then $w_{\infty}$, $w_{\nu}$ for $\nu=0, \cdots, 12$ form an
algebraic equation of degree fourteen, which is just the Jacobian
equation of degree fourteen, whose roots are these $w_{\nu}$ and
$w_{\infty}$:
$$w^{14}+a_1 w^{13}+\cdots+a_{13} w+a_{14}=0.$$

  On the other hand, we have
$$\aligned
  &-13 \sqrt{13} ST^{\nu}(z_1) \cdot ST^{\nu}(z_2) \cdot ST^{\nu}(z_3)\\
 =&-r_4 (\zeta^{8 \nu} z_1^3+\zeta^{7 \nu} z_2^3+\zeta^{11 \nu} z_3^3)
   -r_2 (\zeta^{5 \nu} z_4^3+\zeta^{6 \nu} z_5^3+\zeta^{2 \nu} z_6^3)\\
  &-r_3 (\zeta^{12 \nu} z_1^2 z_2+\zeta^{4 \nu} z_2^2 z_3+\zeta^{10 \nu} z_3^2 z_1)
   -r_1 (\zeta^{\nu} z_4^2 z_5+\zeta^{9 \nu} z_5^2 z_6+\zeta^{3 \nu} z_6^2 z_4)\\
  &+2 r_1 (\zeta^{3 \nu} z_1 z_2^2+\zeta^{\nu} z_2 z_3^2+\zeta^{9 \nu} z_3 z_1^2)
   -2 r_3 (\zeta^{10 \nu} z_4 z_5^2+\zeta^{12 \nu} z_5 z_6^2+\zeta^{4 \nu} z_6 z_4^2)\\
  &+2 r_4 (\zeta^{7 \nu} z_1^2 z_4+\zeta^{11 \nu} z_2^2 z_5+\zeta^{8 \nu} z_3^2 z_6)
   -2 r_2 (\zeta^{6 \nu} z_1 z_4^2+\zeta^{2 \nu} z_2 z_5^2+\zeta^{5 \nu} z_3 z_6^2)+\\
  &+r_1 (\zeta^{3 \nu} z_1^2 z_5+\zeta^{\nu} z_2^2 z_6+\zeta^{9 \nu} z_3^2 z_4)
   +r_3 (\zeta^{10 \nu} z_2 z_4^2+\zeta^{12 \nu} z_3 z_5^2+\zeta^{4 \nu} z_1 z_6^2)+\\
  &+r_2 (\zeta^{6 \nu} z_1^2 z_6+\zeta^{2 \nu} z_2^2 z_4+\zeta^{5 \nu} z_3^2 z_5)
   +r_4 (\zeta^{7 \nu} z_3 z_4^2+\zeta^{11 \nu} z_1 z_5^2+\zeta^{8 \nu} z_2 z_6^2)+\\
  &+r_0 z_1 z_2 z_3+r_{\infty} z_4 z_5 z_6+\\
  &-r_4 (\zeta^{11 \nu} z_1 z_2 z_4+\zeta^{8 \nu} z_2 z_3 z_5+\zeta^{7 \nu} z_1 z_3 z_6)+\\
  &+r_2 (\zeta^{2 \nu} z_1 z_4 z_5+\zeta^{5 \nu} z_2 z_5 z_6+\zeta^{6 \nu} z_3 z_4 z_6)+\\
  &-3 r_4 (\zeta^{7 \nu} z_1 z_2 z_5+\zeta^{11 \nu} z_2 z_3 z_6+\zeta^{8 \nu} z_1 z_3 z_4)+\\
  &+3 r_2 (\zeta^{6 \nu} z_2 z_4 z_5+\zeta^{2 \nu} z_3 z_5 z_6+\zeta^{5 \nu} z_1 z_4 z_6)+\\
  &-r_3 (\zeta^{10 \nu} z_1 z_2 z_6+\zeta^{4 \nu} z_1 z_3 z_5+\zeta^{12 \nu} z_2 z_3 z_4)+\\
  &+r_1 (\zeta^{3 \nu} z_3 z_4 z_5+\zeta^{9 \nu} z_2 z_4 z_6+\zeta^{\nu} z_1 z_5 z_6),
\endaligned$$
where
$$r_0=2(\theta_1-\theta_3)-3(\theta_2-\theta_4), \quad
  r_{\infty}=2(\theta_4-\theta_2)-3(\theta_1-\theta_3),$$
$$r_1=\sqrt{-13-2 \sqrt{13}}, \quad r_2=\sqrt{\frac{-13+3 \sqrt{13}}{2}},$$
$$r_3=\sqrt{-13+2 \sqrt{13}}, \quad r_4=\sqrt{\frac{-13-3 \sqrt{13}}{2}}.$$
This leads us to define the following senary cubic forms (cubic forms in six variables):
$$\left\{\aligned
  {\Bbb D}_0 &=z_1 z_2 z_3,\\
  {\Bbb D}_1 &=2 z_2 z_3^2+z_2^2 z_6-z_4^2 z_5+z_1 z_5 z_6,\\
  {\Bbb D}_2 &=-z_6^3+z_2^2 z_4-2 z_2 z_5^2+z_1 z_4 z_5+3 z_3 z_5 z_6,\\
  {\Bbb D}_3 &=2 z_1 z_2^2+z_1^2 z_5-z_4 z_6^2+z_3 z_4 z_5,\\
  {\Bbb D}_4 &=-z_2^2 z_3+z_1 z_6^2-2 z_4^2 z_6-z_1 z_3 z_5,\\
  {\Bbb D}_5 &=-z_4^3+z_3^2 z_5-2 z_3 z_6^2+z_2 z_5 z_6+3 z_1 z_4 z_6,\\
  {\Bbb D}_6 &=-z_5^3+z_1^2 z_6-2 z_1 z_4^2+z_3 z_4 z_6+3 z_2 z_4 z_5,\\
  {\Bbb D}_7 &=-z_2^3+z_3 z_4^2-z_1 z_3 z_6-3 z_1 z_2 z_5+2 z_1^2 z_4,\\
  {\Bbb D}_8 &=-z_1^3+z_2 z_6^2-z_2 z_3 z_5-3 z_1 z_3 z_4+2 z_3^2 z_6,\\
  {\Bbb D}_9 &=2 z_1^2 z_3+z_3^2 z_4-z_5^2 z_6+z_2 z_4 z_6,\\
  {\Bbb D}_{10} &=-z_1 z_3^2+z_2 z_4^2-2 z_4 z_5^2-z_1 z_2 z_6,\\
  {\Bbb D}_{11} &=-z_3^3+z_1 z_5^2-z_1 z_2 z_4-3 z_2 z_3 z_6+2 z_2^2 z_5,\\
  {\Bbb D}_{12} &=-z_1^2 z_2+z_3 z_5^2-2 z_5 z_6^2-z_2 z_3 z_4,\\
  {\Bbb D}_{\infty}&=z_4 z_5 z_6.
\endaligned\right.\eqno{(3.10)}$$
Then
$$\aligned
  &-13 \sqrt{13} ST^{\nu}({\Bbb D}_0)\\
 =&r_0 {\Bbb D}_0+r_1 \zeta^{\nu} {\Bbb D}_1+r_2 \zeta^{2 \nu} {\Bbb D}_2+
   r_1 \zeta^{3 \nu} {\Bbb D}_3+r_3 \zeta^{4 \nu} {\Bbb D}_4+
   r_2 \zeta^{5 \nu} {\Bbb D}_5+r_2 \zeta^{6 \nu} {\Bbb D}_6+\\
  &+r_4 \zeta^{7 \nu} {\Bbb D}_7+r_4 \zeta^{8 \nu} {\Bbb D}_8
   +r_1 \zeta^{9 \nu} {\Bbb D}_9+r_3 \zeta^{10 \nu} {\Bbb D}_{10}+r_4 \zeta^{11 \nu} {\Bbb D}_{11}+
   r_3 \zeta^{12 \nu} {\Bbb D}_{12}+r_{\infty} {\Bbb D}_{\infty}.
\endaligned$$
$$\aligned
  &-13 \sqrt{13} ST^{\nu}({\Bbb D}_{\infty})\\
 =&r_{\infty} {\Bbb D}_0-r_3 \zeta^{\nu} {\Bbb D}_1-r_4 \zeta^{2 \nu} {\Bbb D}_2-r_3 \zeta^{3 \nu} {\Bbb D}_3
  +r_1 \zeta^{4 \nu} {\Bbb D}_4-r_4 \zeta^{5 \nu} {\Bbb D}_5-r_4 \zeta^{6 \nu} {\Bbb D}_6+\\
  &+r_2 \zeta^{7 \nu} {\Bbb D}_7+r_2 \zeta^{8 \nu} {\Bbb D}_8-r_3 \zeta^{9 \nu} {\Bbb D}_9+r_1 \zeta^{10 \nu}
   {\Bbb D}_{10}+r_2 \zeta^{11 \nu} {\Bbb D}_{11}+r_1 \zeta^{12 \nu} {\Bbb D}_{12}-r_0 {\Bbb D}_{\infty}.
\endaligned$$

  Let
$$\delta_{\infty}(z_1, z_2, z_3, z_4, z_5, z_6)=13^2 (z_1^2 z_2^2 z_3^2+z_4^2 z_5^2 z_6^2)\eqno{(3.11)}$$
and
$$\delta_{\nu}(z_1, z_2, z_3, z_4, z_5, z_6)=\delta_{\infty}(ST^{\nu}(z_1, z_2, z_3, z_4, z_5, z_6))\eqno{(3.12)}$$
for $\nu=0, 1, \cdots, 12$. Then
$$\delta_{\nu}=13^2 ST^{\nu}({\Bbb G}_0)=-13 {\Bbb G}_0+\zeta^{\nu} {\Bbb G}_1
  +\zeta^{2 \nu} {\Bbb G}_2+\cdots+\zeta^{12 \nu} {\Bbb G}_{12},\eqno{(3.13)}$$
where the senary sextic forms (i.e., sextic forms in six
variables) are given as follows:
$$\left\{\aligned
  {\Bbb G}_0 &={\Bbb D}_0^2+{\Bbb D}_{\infty}^2,\\
  {\Bbb G}_1 &=-{\Bbb D}_7^2+2 {\Bbb D}_0 {\Bbb D}_1+10 {\Bbb D}_{\infty} {\Bbb D}_1
               +2 {\Bbb D}_2 {\Bbb D}_{12}-2 {\Bbb D}_3 {\Bbb D}_{11}-4 {\Bbb D}_4 {\Bbb D}_{10}
               -2 {\Bbb D}_9 {\Bbb D}_5,\\
  {\Bbb G}_2 &=-2 {\Bbb D}_1^2-4 {\Bbb D}_0 {\Bbb D}_2+6 {\Bbb D}_{\infty} {\Bbb D}_2
               -2 {\Bbb D}_4 {\Bbb D}_{11}+2 {\Bbb D}_5 {\Bbb D}_{10}-2 {\Bbb D}_6 {\Bbb D}_9
               -2 {\Bbb D}_7 {\Bbb D}_8,\\
  {\Bbb G}_3 &=-{\Bbb D}_8^2+2 {\Bbb D}_0 {\Bbb D}_3+10 {\Bbb D}_{\infty} {\Bbb D}_3
               +2 {\Bbb D}_6 {\Bbb D}_{10}-2 {\Bbb D}_9 {\Bbb D}_7-4 {\Bbb D}_{12} {\Bbb D}_4
               -2 {\Bbb D}_1 {\Bbb D}_2,\\
  {\Bbb G}_4 &=-{\Bbb D}_2^2+10 {\Bbb D}_0 {\Bbb D}_4-2 {\Bbb D}_{\infty} {\Bbb D}_4
               +2 {\Bbb D}_5 {\Bbb D}_{12}-2 {\Bbb D}_9 {\Bbb D}_8-4 {\Bbb D}_1 {\Bbb D}_3
               -2 {\Bbb D}_{10} {\Bbb D}_7,\\
  {\Bbb G}_5 &=-2 {\Bbb D}_9^2-4 {\Bbb D}_0 {\Bbb D}_5+6 {\Bbb D}_{\infty} {\Bbb D}_5
               -2 {\Bbb D}_{10} {\Bbb D}_8+2 {\Bbb D}_6 {\Bbb D}_{12}-2 {\Bbb D}_2 {\Bbb D}_3
               -2 {\Bbb D}_{11} {\Bbb D}_7,\\
  {\Bbb G}_6 &=-2 {\Bbb D}_3^2-4 {\Bbb D}_0 {\Bbb D}_6+6 {\Bbb D}_{\infty} {\Bbb D}_6
               -2 {\Bbb D}_{12} {\Bbb D}_7+2 {\Bbb D}_2 {\Bbb D}_4-2 {\Bbb D}_5 {\Bbb D}_1
               -2 {\Bbb D}_8 {\Bbb D}_{11},\\
  {\Bbb G}_7 &=-2 {\Bbb D}_{10}^2+6 {\Bbb D}_0 {\Bbb D}_7+4 {\Bbb D}_{\infty} {\Bbb D}_7
               -2 {\Bbb D}_1 {\Bbb D}_6-2 {\Bbb D}_2 {\Bbb D}_5-2 {\Bbb D}_8 {\Bbb D}_{12}
               -2 {\Bbb D}_9 {\Bbb D}_{11},\\
  {\Bbb G}_8 &=-2 {\Bbb D}_4^2+6 {\Bbb D}_0 {\Bbb D}_8+4 {\Bbb D}_{\infty} {\Bbb D}_8
               -2 {\Bbb D}_3 {\Bbb D}_5-2 {\Bbb D}_6 {\Bbb D}_2-2 {\Bbb D}_{11} {\Bbb D}_{10}
               -2 {\Bbb D}_1 {\Bbb D}_7,\\
  {\Bbb G}_9 &=-{\Bbb D}_{11}^2+2 {\Bbb D}_0 {\Bbb D}_9+10 {\Bbb D}_{\infty} {\Bbb D}_9
               +2 {\Bbb D}_5 {\Bbb D}_4-2 {\Bbb D}_1 {\Bbb D}_8-4 {\Bbb D}_{10} {\Bbb D}_{12}
               -2 {\Bbb D}_3 {\Bbb D}_6,\\
  {\Bbb G}_{10} &=-{\Bbb D}_5^2+10 {\Bbb D}_0 {\Bbb D}_{10}-2 {\Bbb D}_{\infty} {\Bbb D}_{10}
               +2 {\Bbb D}_6 {\Bbb D}_4-2 {\Bbb D}_3 {\Bbb D}_7-4 {\Bbb D}_9 {\Bbb D}_1
               -2 {\Bbb D}_{12} {\Bbb D}_{11},\\
  {\Bbb G}_{11} &=-2 {\Bbb D}_{12}^2+6 {\Bbb D}_0 {\Bbb D}_{11}+4 {\Bbb D}_{\infty} {\Bbb D}_{11}
               -2 {\Bbb D}_9 {\Bbb D}_2-2 {\Bbb D}_5 {\Bbb D}_6-2 {\Bbb D}_7 {\Bbb D}_4
               -2 {\Bbb D}_3 {\Bbb D}_8,\\
  {\Bbb G}_{12} &=-{\Bbb D}_6^2+10 {\Bbb D}_0 {\Bbb D}_{12}-2 {\Bbb D}_{\infty} {\Bbb D}_{12}
               +2 {\Bbb D}_2 {\Bbb D}_{10}-2 {\Bbb D}_1 {\Bbb D}_{11}-4 {\Bbb D}_3 {\Bbb D}_9
               -2 {\Bbb D}_4 {\Bbb D}_8.
\endaligned\right.\eqno{(3.14)}$$
We have that ${\Bbb G}_0$ is invariant under the action of
$\langle H, T \rangle$, a maximal subgroup of order $78$ of $G$
with index $14$. Note that $\delta_{\infty}$, $\delta_{\nu}$ for
$\nu=0, \cdots, 12$ form an algebraic equation of degree fourteen.
However, we have $\delta_{\infty}+\sum_{\nu=0}^{12}
\delta_{\nu}=0$. Hence, it is not the Jacobian equation of degree
fourteen! We call it the exotic equation of degree fourteen.

  Recall that the theta functions with characteristic
$\left[\matrix \epsilon\\ \epsilon^{\prime} \endmatrix\right] \in
{\Bbb R}^2$ is defined by the following series which converges
uniformly and absolutely on compact subsets of ${\Bbb C} \times
{\Bbb H}$ (see \cite{FK}, p. 73):
$$\theta \left[\matrix \epsilon\\ \epsilon^{\prime} \endmatrix\right] (z, \tau)
 =\sum_{n \in {\Bbb Z}} \exp \left\{2 \pi i \left[\frac{1}{2}
  \left(n+\frac{\epsilon}{2}\right)^2 \tau+\left(n+\frac{\epsilon}{2}\right)
  \left(z+\frac{\epsilon^{\prime}}{2}\right)\right]\right\}.$$
The modified theta constants (see \cite{FK}, p. 215)
$\varphi_l(\tau):=\theta [\chi_l](0, k \tau)$,
where the characteristic $\chi_l=\left[\matrix \frac{2l+1}{k}\\ 1
\endmatrix\right]$, $l=0, \cdots, \frac{k-3}{2}$, for odd $k$ and
$\chi_l=\left[\matrix \frac{2l}{k}\\ 0 \endmatrix\right]$, $l=0,
\cdots, \frac{k}{2}$, for even $k$. We have the following:

{\bf Proposition 3.1.} (see \cite{FK}, p. 236). {\it For each odd
integer $k \geq 5$, the map
$\Phi: \tau \mapsto (\varphi_0(\tau), \varphi_1(\tau), \cdots,
 \varphi_{\frac{k-5}{2}}(\tau), \varphi_{\frac{k-3}{2}}(\tau))$
from ${\Bbb H} \cup {\Bbb Q} \cup \{ \infty \}$ to ${\Bbb
C}^{\frac{k-1}{2}}$, defines a holomorphic mapping from
$\overline{{\Bbb H}/\Gamma(k)}$ into ${\Bbb C} {\Bbb
P}^{\frac{k-3}{2}}$.}

  In our case, the map
$\Phi: \tau \mapsto (\varphi_0(\tau), \varphi_1(\tau), \varphi_2(\tau),
 \varphi_3(\tau), \varphi_4(\tau), \varphi_5(\tau))$
gives a holomorphic mapping from the modular curve
$X(13)=\overline{{\Bbb H}/\Gamma(13)}$ into ${\Bbb C} {\Bbb P}^5$,
which corresponds to our six-dimensional representation, i.e., up
to the constants, $z_1, \cdots, z_6$ are just modular forms
$\varphi_0(\tau), \cdots, \varphi_5(\tau)$. Let
$$\left\{\aligned
  a_1(z) &:=e^{-\frac{11 \pi i}{26}} \theta \left[\matrix \frac{11}{13}\\ 1 \endmatrix\right](0, 13z)
           =q^{\frac{121}{104}} \sum_{n \in {\Bbb Z}} (-1)^n q^{\frac{1}{2}(13n^2+11n)},\\
  a_2(z) &:=e^{-\frac{7 \pi i}{26}} \theta \left[\matrix \frac{7}{13}\\ 1 \endmatrix\right](0, 13z)
           =q^{\frac{49}{104}} \sum_{n \in {\Bbb Z}} (-1)^n q^{\frac{1}{2}(13n^2+7n)},\\
  a_3(z) &:=e^{-\frac{5 \pi i}{26}} \theta \left[\matrix \frac{5}{13}\\ 1 \endmatrix\right](0, 13z)
           =q^{\frac{25}{104}} \sum_{n \in {\Bbb Z}} (-1)^n q^{\frac{1}{2}(13n^2+5n)},\\
  a_4(z) &:=-e^{-\frac{3 \pi i}{26}} \theta \left[\matrix \frac{3}{13}\\ 1 \endmatrix\right](0, 13z)
           =-q^{\frac{9}{104}} \sum_{n \in {\Bbb Z}} (-1)^n q^{\frac{1}{2}(13n^2+3n)},\\
  a_5(z) &:=e^{-\frac{9 \pi i}{26}} \theta \left[\matrix \frac{9}{13}\\ 1 \endmatrix\right](0, 13z)
           =q^{\frac{81}{104}} \sum_{n \in {\Bbb Z}} (-1)^n q^{\frac{1}{2}(13n^2+9n)},\\
  a_6(z) &:=e^{-\frac{\pi i}{26}} \theta \left[\matrix \frac{1}{13}\\ 1 \endmatrix\right](0, 13z)
           =q^{\frac{1}{104}} \sum_{n \in {\Bbb Z}} (-1)^n q^{\frac{1}{2}(13n^2+n)}
\endaligned\right.\eqno{(3.15)}$$
be theta constants of order $13$ and ${\bold A}(z):=(a_1(z), a_2(z), a_3(z), a_4(z), a_5(z), a_6(z))^{T}$.
The significance of our six dimensional representation of $\text{PSL}(2, 13)$ comes from the following:

{\bf Proposition 3.2} (see \cite{Y2}, Proposition 2.5). {\it If $z \in {\Bbb H}$, then the following
relations hold:
$${\bold A}(z+1)=e^{-\frac{3 \pi i}{4}} T {\bold A}(z), \quad
  {\bold A}\left(-\frac{1}{z}\right)=e^{\frac{\pi i}{4}} \sqrt{z} S {\bold A}(z),\eqno{(3.16)}$$
where $S$ and $T$ are given in (3.1) and (3.2), and $0<\text{arg} \sqrt{z} \leq \pi/2$.}

  Recall that the principal congruence subgroup of level $13$ is the normal subgroup $\Gamma(13)$ of
$\Gamma=\text{PSL}(2, {\Bbb Z})$ defined by the exact sequence
$1 \rightarrow \Gamma(13) \rightarrow \Gamma(1) @> f >> G \rightarrow 1$
where $f(\gamma) \equiv \gamma$ (mod $13$) for $\gamma \in \Gamma=\Gamma(1)$. There is a representation
$\rho: \Gamma \rightarrow \text{PGL}(6, {\Bbb C})$ with kernel $\Gamma(13)$ defined as follows: if
$t=\left(\matrix 1 & 1\\ 0 & 1 \endmatrix\right)$ and $s=\left(\matrix 0 & -1\\ 1 & 0 \endmatrix\right)$,
then $\rho(t)=T$ and $\rho(s)=S$. To see that such a representation exists, note that $\Gamma$ is defined by the
presentation $\langle s, t; s^2=(st)^3=1 \rangle$ satisfied by $s$ and $t$ and we have proved that $S$ and $T$
satisfy these relations. Moreover, we have proved that $G$ is defined by the presentation
$\langle S, T; S^2=T^{13}=(ST)^3=1 \rangle$. Let $p=s t^{-1} s$ and $q=s t^3$. Then
$$h:=q^5 p^2 \cdot p^2 q^6 p^8 \cdot q^5 p^2 \cdot p^3 q
    =\left(\matrix
     4, 428, 249 & -10, 547, 030\\
    -11, 594, 791 & 27, 616, 019
     \endmatrix\right)$$
satisfies that $\rho(h)=H$.

  Let $x_i(z)=\eta(z) a_i(z)$, $y_i(z)=\eta^3(z) a_i(z)$ $(1 \leq i \leq 6)$, $X(z)=(x_1(z), \cdots, x_6(z))^{T}$
and $Y(z)=(y_1(z), \cdots, y_6(z))^{T}$. Then $X(z)=\eta(z) {\bold A}(z)$ and $Y(z)=\eta^3(z) {\bold A}(z)$.
Recall that $\eta(z)$ satisfies the following transformation formulas $\eta(z+1)=e^{\frac{\pi i}{12}} \eta(z)$ and
$\eta\left(-\frac{1}{z}\right)=e^{-\frac{\pi i}{4}} \sqrt{z} \eta(z)$. By Proposition 3.2, we have
$X(z+1)=e^{-\frac{2 \pi i}{3}} \rho(t) X(z)$, $X\left(-\frac{1}{z}\right)=z \rho(s) X(z)$,
$Y(z+1)=e^{-\frac{\pi i}{2}} \rho(t) Y(z)$ and $Y\left(-\frac{1}{z}\right)=e^{-\frac{\pi i}{2}} z^2 \rho(s) Y(z)$.
Define $j(\gamma, z):=cz+d$ if $z \in {\Bbb H}$ and $\gamma=\left(\matrix a & b\\ c & d \endmatrix\right) \in \Gamma(1)$.
Hence, $X(\gamma(z))=u(\gamma) j(\gamma, z) \rho(\gamma) X(z)$ and $Y(\gamma(z))=v(\gamma) j(\gamma, z)^2 \rho(\gamma) Y(z)$
for $\gamma \in \Gamma(1)$, where $u(\gamma)=1, \omega$ or $\omega^2$ with $\omega=e^{\frac{2 \pi i}{3}}$ and
$v(\gamma)=\pm 1$ or $\pm i$. Since $\Gamma(13)=\text{ker}$ $\rho$, we have
$X(\gamma(z))=u(\gamma) j(\gamma, z) X(z)$ and $Y(\gamma(z))=v(\gamma) j(\gamma, z)^2 Y(z)$ for $\gamma \in \Gamma(13)$.
This means that the functions $x_1(z)$, $\cdots$, $x_6(z)$ are modular forms of weight one for
$\Gamma(13)$ with the same multiplier $u(\gamma)=1, \omega$ or $\omega^2$ and $y_1(z)$, $\cdots$, $y_6(z)$ are modular
forms of weight two for $\Gamma(13)$ with the same multiplier $v(\gamma)=\pm 1 $ or $\pm i$.

  From now on, we will use the following abbreviation ${\Bbb A}_j={\Bbb A_j}(a_1(z), \cdots, a_6(z))$ $(0 \leq j \leq
6)$, ${\Bbb D}_j={\Bbb D_j}(a_1(z), \cdots, a_6(z))$ $(j=0,1, \cdots, 12, \infty)$ and ${\Bbb G}_j={\Bbb G_j}(a_1(z),
\cdots, a_6(z))$ $(0 \leq j \leq 12)$. We have
$${\Bbb A}_0=q^{\frac{1}{4}} (1+O(q)),\quad
  {\Bbb A}_1=q^{\frac{17}{52}} (2+O(q)),\quad
  {\Bbb A}_2=q^{\frac{29}{52}} (2+O(q)),\quad
  {\Bbb A}_3=q^{\frac{49}{52}} (1+O(q)),$$
$${\Bbb A}_4=q^{\frac{25}{52}} (-1+O(q)),\quad
  {\Bbb A}_5=q^{\frac{9}{52}} (-1+O(q)),\quad
  {\Bbb A}_6=q^{\frac{1}{52}} (-1+O(q)),$$
and
$$\left\{\aligned
  {\Bbb D}_0 &=q^{\frac{15}{8}} (1+O(q)),\\
  {\Bbb D}_{\infty} &=q^{\frac{7}{8}} (-1+O(q)),\\
  {\Bbb D}_1 &=q^{\frac{99}{104}} (2+O(q)),\\
  {\Bbb D}_2 &=q^{\frac{3}{104}} (-1+O(q)),\\
  {\Bbb D}_3 &=q^{\frac{11}{104}} (1+O(q)),\\
  {\Bbb D}_4 &=q^{\frac{19}{104}} (-2+O(q)),\\
  {\Bbb D}_5 &=q^{\frac{27}{104}} (-1+O(q)),
\endaligned\right. \quad \quad
  \left\{\aligned
  {\Bbb D}_6 &=q^{\frac{35}{104}} (-1+O(q)),\\
  {\Bbb D}_7 &=q^{\frac{43}{104}} (1+O(q)),\\
  {\Bbb D}_8 &=q^{\frac{51}{104}} (3+O(q)),\\
  {\Bbb D}_9 &=q^{\frac{59}{104}} (-2+O(q)),\\
  {\Bbb D}_{10} &=q^{\frac{67}{104}} (1+O(q)),\\
  {\Bbb D}_{11} &=q^{\frac{75}{104}} (-4+O(q)),\\
  {\Bbb D}_{12} &=q^{\frac{83}{104}} (-1+O(q)).
\endaligned\right.$$
Hence,
$$\left\{\aligned
  {\Bbb G}_0 &=q^{\frac{7}{4}} (1+O(q)),\\
  {\Bbb G}_1 &=q^{\frac{43}{52}} (13+O(q)),\\
  {\Bbb G}_2 &=q^{\frac{47}{52}} (-22+O(q)),\\
  {\Bbb G}_3 &=q^{\frac{51}{52}} (-21+O(q)),\\
  {\Bbb G}_4 &=q^{\frac{3}{52}} (-1+O(q)),\\
  {\Bbb G}_5 &=q^{\frac{7}{52}} (2+O(q)),\\
  {\Bbb G}_6 &=q^{\frac{11}{52}} (2+O(q)),
\endaligned\right. \quad \quad
  \left\{\aligned
  {\Bbb G}_7 &=q^{\frac{15}{52}} (-2+O(q)),\\
  {\Bbb G}_8 &=q^{\frac{19}{52}} (-8+O(q)),\\
  {\Bbb G}_9 &=q^{\frac{23}{52}} (6+O(q)),\\
  {\Bbb G}_{10} &=q^{\frac{27}{52}} (1+O(q)),\\
  {\Bbb G}_{11} &=q^{\frac{31}{52}} (-8+O(q)),\\
  {\Bbb G}_{12} &=q^{\frac{35}{52}} (17+O(q)).
\endaligned\right.$$
Note that
$$\aligned
  w_{\nu} &=({\Bbb A}_0+\zeta^{\nu} {\Bbb A}_1+\zeta^{4 \nu} {\Bbb A}_2+\zeta^{9 \nu} {\Bbb A}_3+
           +\zeta^{3 \nu} {\Bbb A}_4+\zeta^{12 \nu} {\Bbb A}_5+\zeta^{10 \nu} {\Bbb A}_6)^2\\
          &={\Bbb A}_0^2+2 ({\Bbb A}_1 {\Bbb A}_5+{\Bbb A}_2 {\Bbb A}_3+{\Bbb A}_4 {\Bbb A}_6)+\\
          &+2 \zeta^{\nu} ({\Bbb A}_0 {\Bbb A}_1+{\Bbb A}_2 {\Bbb A}_6)
           +2 \zeta^{3 \nu} ({\Bbb A}_0 {\Bbb A}_4+{\Bbb A}_2 {\Bbb A}_5)
           +2 \zeta^{9 \nu} ({\Bbb A}_0 {\Bbb A}_3+{\Bbb A}_5 {\Bbb A}_6)+\\
          &+2 \zeta^{12 \nu} ({\Bbb A}_0 {\Bbb A}_5+{\Bbb A}_3 {\Bbb A}_4)
           +2 \zeta^{10 \nu} ({\Bbb A}_0 {\Bbb A}_6+{\Bbb A}_1 {\Bbb A}_3)
           +2 \zeta^{4 \nu} ({\Bbb A}_0 {\Bbb A}_2+{\Bbb A}_1 {\Bbb A}_4)+\\
          &+\zeta^{2 \nu} ({\Bbb A}_1^2+2 {\Bbb A}_4 {\Bbb A}_5)
           +\zeta^{5 \nu} ({\Bbb A}_3^2+2 {\Bbb A}_1 {\Bbb A}_2)
           +\zeta^{6 \nu} ({\Bbb A}_4^2+2 {\Bbb A}_3 {\Bbb A}_6)+\\
          &+\zeta^{11 \nu} ({\Bbb A}_5^2+2 {\Bbb A}_1 {\Bbb A}_6)
           +\zeta^{8 \nu} ({\Bbb A}_2^2+2 {\Bbb A}_3 {\Bbb A}_5)
           +\zeta^{7 \nu} ({\Bbb A}_6^2+2 {\Bbb A}_4 {\Bbb A}_2),
\endaligned$$
where
$$\left\{\aligned
  {\Bbb A}_0^2+2 ({\Bbb A}_1 {\Bbb A}_5+{\Bbb A}_2 {\Bbb A}_3+{\Bbb A}_4 {\Bbb A}_6)&=q^{\frac{1}{2}} (-1+O(q)),\\
  {\Bbb A}_0 {\Bbb A}_1+{\Bbb A}_2 {\Bbb A}_6 &=q^{\frac{41}{26}} (-3+O(q)),\\
  {\Bbb A}_0 {\Bbb A}_4+{\Bbb A}_2 {\Bbb A}_5 &=q^{\frac{19}{26}} (-3+O(q)),\\
  {\Bbb A}_0 {\Bbb A}_3+{\Bbb A}_5 {\Bbb A}_6 &=q^{\frac{5}{26}} (1+O(q)),\\
  {\Bbb A}_0 {\Bbb A}_5+{\Bbb A}_3 {\Bbb A}_4 &=q^{\frac{11}{26}} (-1+O(q)),\\
  {\Bbb A}_0 {\Bbb A}_6+{\Bbb A}_1 {\Bbb A}_3 &=q^{\frac{7}{26}} (-1+O(q)),\\
  {\Bbb A}_0 {\Bbb A}_2+{\Bbb A}_1 {\Bbb A}_4 &=q^{\frac{47}{26}} (-1+O(q)),\\
  {\Bbb A}_1^2+2 {\Bbb A}_4 {\Bbb A}_5 &=q^{\frac{17}{26}} (6+O(q)),\\
  {\Bbb A}_3^2+2 {\Bbb A}_1 {\Bbb A}_2 &=q^{\frac{23}{26}} (8+O(q)),\\
  {\Bbb A}_4^2+2 {\Bbb A}_3 {\Bbb A}_6 &=q^{\frac{25}{26}} (-1+O(q)),\\
  {\Bbb A}_5^2+2 {\Bbb A}_1 {\Bbb A}_6 &=q^{\frac{9}{26}} (-3+O(q)),\\
  {\Bbb A}_2^2+2 {\Bbb A}_3 {\Bbb A}_5 &=q^{\frac{29}{26}} (2+O(q)),\\
  {\Bbb A}_6^2+2 {\Bbb A}_4 {\Bbb A}_2 &=q^{\frac{1}{26}} (1+O(q)).
\endaligned\right.$$

{\it Proof of Theorem 1.1}. Let
$$\Phi_{20}=w_0^5+w_1^5+\cdots+w_{12}^5+w_{\infty}^5.$$
As a polynomial in six variables, $\Phi_{20}(z_1, z_2, z_3, z_4, z_5, z_6)$ is
a $G$-invariant polynomial. Moreover, for $\gamma \in \Gamma(1)$,
$$\aligned
  &\Phi_{20}(Y(\gamma(z))^{T})=\Phi_{20}(v(\gamma) j(\gamma, z)^2 (\rho(\gamma) Y(z))^{T})\\
 =&v(\gamma)^{20} j(\gamma, z)^{40} \Phi_{20}((\rho(\gamma) Y(z))^{T})=j(\gamma, z)^{40} \Phi_{20}((\rho(\gamma) Y(z))^{T}).
\endaligned$$
Note that $\rho(\gamma) \in \langle \rho(s), \rho(t) \rangle=G$ and $\Phi_{20}$ is a $G$-invariant polynomial, we have
$$\Phi_{20}(Y(\gamma(z))^{T})=j(\gamma, z)^{40} \Phi_{20}(Y(z)^{T}), \quad \text{for $\gamma \in \Gamma(1)$}.$$
This implies that $\Phi_{20}(y_1(z), \cdots, y_6(z))$ is a modular form of weight $40$ for the full modular group $\Gamma(1)$.
Moreover, we will show that it is a cusp form. In fact,
$$\aligned
  &\Phi_{20}(a_1(z), \cdots, a_6(z))=13^5 q^{\frac{5}{2}} (1+O(q))^5+\\
  &+\sum_{\nu=0}^{12} [q^{\frac{1}{2}} (-1+O(q))+\\
  &+2 \zeta^{\nu} q^{\frac{41}{26}} (-3+O(q))+2 \zeta^{3 \nu} q^{\frac{19}{26}} (-3+O(q))
   +2 \zeta^{9 \nu} q^{\frac{5}{26}} (1+O(q))+\\
  &+2 \zeta^{12 \nu} q^{\frac{11}{26}} (-1+O(q))+2 \zeta^{10 \nu} q^{\frac{7}{26}} (-1+O(q))
   +2 \zeta^{4 \nu} q^{\frac{47}{26}} (-1+O(q))+\\
  &+\zeta^{2 \nu} q^{\frac{17}{26}} (6+O(q))+\zeta^{5 \nu} q^{\frac{23}{26}} (8+O(q))
   +\zeta^{6 \nu} q^{\frac{25}{26}} (-1+O(q))+\\
  &+\zeta^{11 \nu} q^{\frac{9}{26}} (-3+O(q))+\zeta^{8 \nu} q^{\frac{29}{26}} (2+O(q))
   +\zeta^{7 \nu} q^{\frac{1}{26}} (1+O(q))]^5.
\endaligned$$
We will calculate the $q^{\frac{1}{2}}$-term which is the lowest degree.
For the partition $13=4 \cdot 1+9$, the corresponding term is
$$\left(\matrix 5\\ 4, 1 \endmatrix\right) (\zeta^{7 \nu} q^{\frac{1}{26}})^4 \cdot (-3)
  \zeta^{11 \nu} q^{\frac{9}{26}}=-15 q^{\frac{1}{2}}.$$
For the partition $13=3 \cdot 1+2 \cdot 5$, the corresponding term is
$$\left(\matrix 5\\ 3, 2 \endmatrix\right) (\zeta^{7 \nu} q^{\frac{1}{26}})^3 \cdot
  (2 \zeta^{9 \nu} q^{\frac{5}{26}})^2=40 q^{\frac{1}{2}}.$$
Hence, for $\Phi_{20}(y_1(z), \cdots, y_6(z))$ which is a modular form for $\Gamma(1)$ with weight $40$,
the lowest degree term is given by
$$(-15+40) q^{\frac{1}{2}} \cdot q^{\frac{3}{24} \cdot 20}=25 q^3.$$
Thus,
$$\Phi_{20}(y_1(z), \cdots, y_6(z))=q^3 (13 \cdot 25+O(q)).$$
The leading term of $\Phi_{20}(y_1(z), \cdots, y_6(z))$ together with its weight $40$ suffice to identify
this modular form with $\Phi_{20}(y_1(z), \cdots, y_6(z))=13 \cdot 25 \Delta(z)^3 E_4(z)$. Consequently,
$$\Phi_{20}(x_1(z), \cdots, x_6(z))=13 \cdot 25 \Delta(z)^3 E_4(z)/\eta(z)^{40}=13 \cdot 25 \eta(z)^8 \Delta(z) E_4(z).$$

  Let
$$\Phi_{18}=\delta_0^3+\delta_1^3+\cdots+\delta_{12}^3+\delta_{\infty}^3.$$
As a polynomial in six variables, $\Phi_{18}(z_1, z_2, z_3, z_4, z_5, z_6)$ is
a $G$-invariant polynomial. Moreover, for $\gamma \in \Gamma(1)$,
$$\aligned
  &\Phi_{18}(X(\gamma(z))^{T})=\Phi_{18}(u(\gamma) j(\gamma, z) (\rho(\gamma) X(z))^{T})\\
 =&u(\gamma)^{18} j(\gamma, z)^{18} \Phi_{18}((\rho(\gamma) X(z))^{T})=j(\gamma, z)^{18} \Phi_{18}((\rho(\gamma) X(z))^{T}).
\endaligned$$
Note that $\rho(\gamma) \in \langle \rho(s), \rho(t) \rangle=G$ and $\Phi_{18}$ is a $G$-invariant polynomial, we have
$$\Phi_{18}(X(\gamma(z))^{T})=j(\gamma, z)^{18} \Phi_{18}(X(z)^{T}), \quad \text{for $\gamma \in \Gamma(1)$}.$$
This implies that $\Phi_{18}(x_1(z), \cdots, x_6(z))$ is a modular form of weight $18$ for the full modular group $\Gamma(1)$.
Moreover, we will show that it is a cusp form. In fact,
$$\aligned
  &\Phi_{18}(a_1(z), \cdots, a_6(z))=13^6 q^{\frac{21}{4}} (1+O(q))^3+\\
  &+\sum_{\nu=0}^{12} [-13 q^{\frac{7}{4}} (1+O(q))+\\
  &+\zeta^{\nu} q^{\frac{43}{52}} (13+O(q))+\zeta^{2 \nu} q^{\frac{47}{52}} (-22+O(q))
   +\zeta^{3 \nu} q^{\frac{51}{52}} (-21+O(q))+\\
  &+\zeta^{4 \nu} q^{\frac{3}{52}} (-1+O(q))+\zeta^{5 \nu} q^{\frac{7}{52}} (2+O(q))
   +\zeta^{6 \nu} q^{\frac{11}{52}} (2+O(q))+\\
  &+\zeta^{7 \nu} q^{\frac{15}{52}} (-2+O(q))+\zeta^{8 \nu} q^{\frac{19}{52}} (-8+O(q))
   +\zeta^{9 \nu} q^{\frac{23}{52}} (6+O(q))+\\
  &+\zeta^{10 \nu} q^{\frac{27}{52}} (1+O(q))+\zeta^{11 \nu} q^{\frac{31}{52}} (-8+O(q))
   +\zeta^{12 \nu} q^{\frac{35}{52}} (17+O(q))]^3.
\endaligned$$
We will calculate the $q^{\frac{1}{4}}$-term which is the lowest degree.
For the partition $13=2 \cdot 3+7$, the corresponding term is
$$\left(\matrix 3\\ 2, 1 \endmatrix\right) (\zeta^{4 \nu} q^{\frac{3}{52}} \cdot(-1))^2 \cdot
  \zeta^{5 \nu} q^{\frac{7}{52}} \cdot 2=6 q^{\frac{1}{4}}.$$
Hence, for $\Phi_{18}(x_1(z), \cdots, x_6(z))$ which is a modular form for $\Gamma(1)$ with weight $18$,
the lowest degree term is given by $6 q^{\frac{1}{4}} \cdot q^{\frac{18}{24}}=6 q$.
Thus,
$$\Phi_{18}(x_1(z), \cdots, x_6(z))=q (13 \cdot 6+O(q)).$$
The leading term of $\Phi_{18}(x_1(z), \cdots, x_6(z))$ together with its weight $18$ suffice to identify
this modular form with
$$\Phi_{18}(x_1(z), \cdots, x_6(z))=13 \cdot 6 \Delta(z) E_6(z).$$

  Let
$$\Phi_{12}=\delta_0^2+\delta_1^2+\cdots+\delta_{12}^2+\delta_{\infty}^2.$$
As a polynomial in six variables, $\Phi_{12}(z_1, z_2, z_3, z_4, z_5, z_6)$ is
a $G$-invariant polynomial. Similarly as above, we can show that $\Phi_{12}(x_1(z), \cdots, x_6(z))$ is a
modular form of weight $12$ for the full modular group $\Gamma(1)$. Moreover, we will show that it is a
cusp form. In fact,
$$\aligned
  &\Phi_{12}(a_1(z), \cdots, a_6(z))=13^4 q^{\frac{7}{2}} (1+O(q))^2+\\
  &+\sum_{\nu=0}^{12} [-13 q^{\frac{7}{4}} (1+O(q))+\\
  &+\zeta^{\nu} q^{\frac{43}{52}} (13+O(q))+\zeta^{2 \nu} q^{\frac{47}{52}} (-22+O(q))
   +\zeta^{3 \nu} q^{\frac{51}{52}} (-21+O(q))+\\
  &+\zeta^{4 \nu} q^{\frac{3}{52}} (-1+O(q))+\zeta^{5 \nu} q^{\frac{7}{52}} (2+O(q))
   +\zeta^{6 \nu} q^{\frac{11}{52}} (2+O(q))+\\
  &+\zeta^{7 \nu} q^{\frac{15}{52}} (-2+O(q))+\zeta^{8 \nu} q^{\frac{19}{52}} (-8+O(q))
   +\zeta^{9 \nu} q^{\frac{23}{52}} (6+O(q))+\\
  &+\zeta^{10 \nu} q^{\frac{27}{52}} (1+O(q))+\zeta^{11 \nu} q^{\frac{31}{52}} (-8+O(q))
   +\zeta^{12 \nu} q^{\frac{35}{52}} (17+O(q))]^2.
\endaligned$$
We will calculate the $q^{\frac{1}{2}}$-term which is the lowest degree.
For the partition $26=3+23$, the corresponding term is
$$\left(\matrix 2\\ 1, 1 \endmatrix\right) \zeta^{4 \nu} q^{\frac{3}{52}} \cdot (-1) \cdot
  \zeta^{9 \nu} q^{\frac{23}{52}} \cdot 6=-12 q^{\frac{1}{2}}.$$
For the partition $26=7+19$, the corresponding term is
$$\left(\matrix 2\\ 1, 1 \endmatrix\right) \zeta^{5 \nu} q^{\frac{7}{52}} \cdot 2 \cdot
  \zeta^{8 \nu} q^{\frac{19}{52}} \cdot (-8)=-32 q^{\frac{1}{2}}.$$
For the partition $26=11+15$, the corresponding term is
$$\left(\matrix 2\\ 1, 1 \endmatrix\right) \zeta^{6 \nu} q^{\frac{11}{52}} \cdot 2 \cdot
  \zeta^{7 \nu} q^{\frac{15}{52}} \cdot (-2)=-8 q^{\frac{1}{2}}.$$
Hence, for $\Phi_{12}(x_1(z), \cdots, x_6(z))$ which is a modular form for $\Gamma(1)$ with weight $12$,
the lowest degree term is given by $(-12-32-8) q^{\frac{1}{2}} \cdot q^{\frac{12}{24}}=-52 q$.
Thus,
$$\Phi_{12}(x_1(z), \cdots, x_6(z))=q (-13 \cdot 52+O(q)).$$
The leading term of $\Phi_{12}(x_1(z), \cdots, x_6(z))$ together with its weight $12$ suffice to identify
this modular form with
$$\Phi_{12}(x_1(z), \cdots, x_6(z))=-13 \cdot 52 \Delta(z).$$

  Let
$$\Phi_{30}=\delta_0^5+\delta_1^5+\cdots+\delta_{12}^5+\delta_{\infty}^5.$$
As a polynomial in six variables, $\Phi_{30}(z_1, z_2, z_3, z_4, z_5, z_6)$ is
a $G$-invariant polynomial. Similarly as above, we can show that $\Phi_{30}(x_1(z), \cdots, x_6(z))$ is a
modular form of weight $30$ for the full modular group $\Gamma(1)$. Moreover, we will show that it is a
cusp form. In fact,
$$\aligned
  &\Phi_{30}(a_1(z), \cdots, a_6(z))=13^{10} q^{\frac{35}{4}} (1+O(q))^5+\\
  &+\sum_{\nu=0}^{12} [-13 q^{\frac{7}{4}} (1+O(q))+\\
  &+\zeta^{\nu} q^{\frac{43}{52}} (13+O(q))+\zeta^{2 \nu} q^{\frac{47}{52}} (-22+O(q))
   +\zeta^{3 \nu} q^{\frac{51}{52}} (-21+O(q))+\\
  &+\zeta^{4 \nu} q^{\frac{3}{52}} (-1+O(q))+\zeta^{5 \nu} q^{\frac{7}{52}} (2+O(q))
   +\zeta^{6 \nu} q^{\frac{11}{52}} (2+O(q))+\\
  &+\zeta^{7 \nu} q^{\frac{15}{52}} (-2+O(q))+\zeta^{8 \nu} q^{\frac{19}{52}} (-8+O(q))
   +\zeta^{9 \nu} q^{\frac{23}{52}} (6+O(q))+\\
  &+\zeta^{10 \nu} q^{\frac{27}{52}} (1+O(q))+\zeta^{11 \nu} q^{\frac{31}{52}} (-8+O(q))
   +\zeta^{12 \nu} q^{\frac{35}{52}} (17+O(q))]^5.
\endaligned$$
We will calculate the $q^{\frac{3}{4}}$-term which is the lowest degree.
(1) For the partition $39=4 \cdot 3+27$, the corresponding term is
$$\left(\matrix 5\\ 4, 1 \endmatrix\right) (\zeta^{4 \nu} q^{\frac{3}{52}} \cdot (-1))^4 \cdot
  \zeta^{10 \nu} q^{\frac{27}{52}}=5 q^{\frac{3}{4}}.$$
(2) For the partition $39=3 \cdot 3+7+23$, the corresponding term is
$$\left(\matrix 5\\ 3, 1, 1 \endmatrix\right) (\zeta^{4 \nu} q^{\frac{3}{52}} \cdot (-1))^3 \cdot
  \zeta^{5 \nu} q^{\frac{7}{52}} \cdot 2 \cdot \zeta^{9 \nu} q^{\frac{23}{52}} \cdot 6=-240 q^{\frac{3}{4}}.$$
(3) For the partition $39=3 \cdot 3+11+19$, the corresponding term is
$$\left(\matrix 5\\ 3, 1, 1 \endmatrix\right) (\zeta^{4 \nu} q^{\frac{3}{52}} \cdot (-1))^3 \cdot
  \zeta^{6 \nu} q^{\frac{11}{52}} \cdot 2 \cdot \zeta^{8 \nu} q^{\frac{19}{52}} \cdot (-8)=320 q^{\frac{3}{4}}.$$
(4) For the partition $39=3 \cdot 3+2 \cdot 15$, the corresponding term is
$$\left(\matrix 5\\ 3, 2 \endmatrix\right) (\zeta^{4 \nu} q^{\frac{3}{52}} \cdot (-1))^3 \cdot
  (\zeta^{7 \nu} q^{\frac{15}{52}} \cdot (-2))^2=-40 q^{\frac{3}{4}}.$$
(5) For the partition $39=2 \cdot 3+3 \cdot 11$, the corresponding term is
$$\left(\matrix 5\\ 2, 3 \endmatrix\right) (\zeta^{4 \nu} q^{\frac{3}{52}} \cdot (-1))^2 \cdot
  (\zeta^{6 \nu} q^{\frac{11}{52}} \cdot 2)^3=80 q^{\frac{3}{4}}.$$
(6) For the partition $39=2 \cdot 3+2 \cdot 7+19$, the corresponding term is
$$\left(\matrix 5\\ 2, 2, 1 \endmatrix\right) (\zeta^{4 \nu} q^{\frac{3}{52}} \cdot (-1))^2 \cdot
  (\zeta^{5 \nu} q^{\frac{7}{52}} \cdot 2)^2 \cdot \zeta^{8 \nu} q^{\frac{19}{52}} \cdot (-8)=-960 q^{\frac{3}{4}}.$$
(7) For the partition $39=2 \cdot 3+7+11+15$, the corresponding term is
$$\left(\matrix 5\\ 2, 1, 1, 1 \endmatrix\right) (\zeta^{4 \nu} q^{\frac{3}{52}} \cdot (-1))^2 \cdot
  \zeta^{5 \nu} q^{\frac{7}{52}} \cdot 2 \cdot \zeta^{6 \nu} q^{\frac{11}{52}} \cdot 2 \cdot \zeta^{7 \nu}
  q^{\frac{15}{52}} \cdot (-2)=-480 q^{\frac{3}{4}}.$$
(8) For the partition $39=1 \cdot 3+3 \cdot 7+15$, the corresponding term is
$$\left(\matrix 5\\ 1, 3, 1 \endmatrix\right) \zeta^{4 \nu} q^{\frac{3}{52}} \cdot (-1) \cdot
  (\zeta^{5 \nu} q^{\frac{7}{52}} \cdot 2)^3 \cdot \zeta^{7 \nu} q^{\frac{15}{52}} \cdot (-2)=320 q^{\frac{3}{4}}.$$
(9) For the partition $39=1 \cdot 3+2 \cdot 7+2 \cdot 11$, the corresponding term is
$$\left(\matrix 5\\ 1, 2, 2 \endmatrix\right) \zeta^{4 \nu} q^{\frac{3}{52}} \cdot (-1) \cdot
  (\zeta^{5 \nu} q^{\frac{7}{52}} \cdot 2)^2 \cdot (\zeta^{6 \nu} q^{\frac{11}{52}} \cdot 2)^2=-480 q^{\frac{3}{4}}.$$
(10) For the partition $39=4 \cdot 7+11$, the corresponding term is
$$\left(\matrix 5\\ 4, 1 \endmatrix\right) (\zeta^{5 \nu} q^{\frac{7}{52}} \cdot 2)^4 \cdot
  \zeta^{6 \nu} q^{\frac{11}{52}} \cdot 2=160 q^{\frac{3}{4}}.$$
Hence, for $\Phi_{30}(x_1(z), \cdots, x_6(z))$ which is a modular form for $\Gamma(1)$ with weight $30$,
the lowest degree term is given by
$$(5-240+320-40+80-960-480+320-480+160) q^{\frac{3}{4}} \cdot q^{\frac{30}{24}}=-1315 q^2.$$
Thus,
$$\Phi_{30}(x_1(z), \cdots, x_6(z))=q^2 (-13 \cdot 1315+O(q)).$$
The leading term of $\Phi_{30}(x_1(z), \cdots, x_6(z))$ together with its weight $30$ suffice to identify
this modular form with
$$\Phi_{30}(x_1(z), \cdots, x_6(z))=-13 \cdot 1315 \Delta(z)^2 E_6(z).$$

  Up to a constant, we revise the definition of $\Phi_{12}$, $\Phi_{18}$, $\Phi_{20}$ and $\Phi_{30}$:
$$\Phi_{12}=-\frac{1}{13 \cdot 52} \left(\sum_{\nu=0}^{12} \delta_{\nu}^2+\delta_{\infty}^2\right), \quad
  \Phi_{18}=\frac{1}{13 \cdot 6} \left(\sum_{\nu=0}^{12} \delta_{\nu}^3+\delta_{\infty}^3\right),\eqno{(3.17)}$$
$$\Phi_{20}=\frac{1}{13 \cdot 25} \left(\sum_{\nu=0}^{12} w_{\nu}^5+w_{\infty}^5\right), \quad
  \Phi_{30}=-\frac{1}{13 \cdot 1315} \left(\sum_{\nu=0}^{12} \delta_{\nu}^5+\delta_{\infty}^5\right).\eqno{(3.18)}$$
Consequently,
$$\left\{\aligned
  \Phi_{12}(x_1(z), \cdots, x_6(z)) &=\Delta(z),\\
  \Phi_{18}(x_1(z), \cdots, x_6(z)) &=\Delta(z) E_6(z),\\
  \Phi_{20}(x_1(z), \cdots, x_6(z)) &=\eta(z)^8 \Delta(z) E_4(z),\\
  \Phi_{30}(x_1(z), \cdots, x_6(z)) &=\Delta(z)^2 E_6(z).
\endaligned\right.\eqno{(3.19)}$$
From now on, we will use the following abbreviation $\Phi_j=\Phi_j(x_1(z), \cdots, x_6(z))$ for
$j=12, 18, 20$ and $30$. The relations
$$j(z):=\frac{E_4(z)^3}{\Delta(z)}=\frac{\Phi_{20}^3}{\Phi_{12}^5}, \quad
  j(z)-1728=\frac{E_6(z)^2}{\Delta(z)}=\frac{\Phi_{30}^2}{\Phi_{12}^5}=\frac{\Phi_{18}^2}{\Phi_{12}^3}\eqno{(3.20)}$$
give the equations
$$\Phi_{20}^3-\Phi_{30}^2=1728 \Phi_{12}^5, \quad \Phi_{20}^3-\Phi_{12}^2 \Phi_{18}^2=1728 \Phi_{12}^5.\eqno{(3.21)}$$
This completes the proof of Theorem 1.1.

\flushpar $\qquad \qquad \qquad \qquad \qquad \qquad \qquad \qquad
\qquad \qquad \qquad \qquad \qquad \qquad \qquad \qquad \qquad
\qquad \quad \boxed{}$

  Let us recall some facts about exotic spheres (see \cite{Hi2}). A $k$-dimensional compact oriented
differentiable manifold is called a $k$-sphere if it is homeomorphic to the $k$-dimensional standard
sphere. A $k$-sphere not diffeomorphic to the standard $k$-sphere is said to be exotic. The first exotic
sphere was discovered by Milnor in 1956 (see \cite{Mi1} and \cite{Mi4}). Two $k$-spheres are called
equivalent if there exists an orientation preserving diffeomorphism between them. The equivalence classes
of $k$-spheres constitute for $k \geq 5$ a finite abelian group $\Theta_k$ under the connected sum operation.
$\Theta_k$ contains the subgroup $b P_{k+1}$ of those $k$-spheres which bound a parallelizable manifold.
$b P_{4m}$ ($m \geq 2$) is cyclic of order $2^{2m-2}(2^{2m-1}-1)$ numerator $(4 B_m/m)$, where $B_m$ is
the $m$-th Bernoulli number. Let $g_m$ be the Milnor generator of $b P_{4m}$. If a $(4m-1)$-sphere $\Sigma$
bounds a parallelizable manifold $B$ of dimension $4m$, then the signature $\tau(B)$ of the intersection
form of $B$ is divisible by $8$ and $\Sigma=\frac{\tau(B)}{8} g_m$. For $m=2$ we have $b P_8=\Theta_7=
{\Bbb Z} /28 {\Bbb Z}$. All these results are due to Milnor-Kervaire (see \cite{KM}). In particular,
$$\sum_{i=0}^{2m} z_i \overline{z_i}=1, \quad
  z_0^3+z_1^{6k-1}+z_2^2+\cdots+z_{2m}^2=0$$
is a $(4m-1)$-sphere embedded in $S^{4m+1} \subset {\Bbb C}^{2n+1}$ which represents the element
$(-1)^m k \cdot g_m \in b P_{4m}$. For $m=2$ and $k=1, 2, \cdots, 28$ we get the $28$ classes of
$7$-spheres. Theorem 1.1 shows that the higher dimensional liftings of two distinct symmetry groups
and modular interpretations on the equation of $E_8$-singularity give the same Milnor's standard
generator of $\Theta_7$.

  In fact, $G$ is the symmetry group of regular maps $\{ 7, 3 \}_{13}$, $\{ 13, 3 \}_{7}$ and
$\{ 13, 7 \}_{3}$, which are the generalizations of regular polyhedra in ${\Bbb R}^3$ (see \cite{CoM}, p.140).
$$\matrix
 \text{Maps} & \text{Vertices} & \text{Edges} & \text{Faces} & \text{Characteristic} & \text{Group}\\
 \{ 7, 3 \}_{13} & 182 & 273 & 78 & -13 & \text{PSL}(2, 13)\\
 \{ 13, 3 \}_{7} & 182 & 273 & 42 & -49 & \text{PSL}(2, 13)\\
 \{ 13, 7 \}_{3} & 78 & 273 & 42 & -153 & \text{PSL}(2, 13)
\endmatrix$$

\vskip 0.5 cm

\centerline{\bf 4. Modular parametrizations for some exceptional singularities and}
\centerline{\bf Fermat-Catalan conjecture}

\vskip 0.5 cm

{\it Proof of Theorem 1.2 and Corollary 1.3}. Let
$$\Phi_{32}=w_0^8+w_1^8+\cdots+w_{12}^8+w_{\infty}^8.$$
As a polynomial in six variables, $\Phi_{32}(z_1, z_2, z_3, z_4, z_5, z_6)$ is a $G$-invariant polynomial.
Similarly as above, we can show that $\Phi_{32}(y_1(z), \cdots, y_6(z))$ is a modular form of weight $64$
for the full modular group $\Gamma(1)$. Moreover, we will show that it is a cusp form. In fact,
$$\aligned
  &\Phi_{32}(a_1(z), \cdots, a_6(z))=13^8 q^4 (1+O(q))^8+\\
  &+\sum_{\nu=0}^{12} [q^{\frac{1}{2}} (-1+O(q))+\\
  &+2 \zeta^{\nu} q^{\frac{41}{26}} (-3+O(q))+2 \zeta^{3 \nu} q^{\frac{19}{26}} (-3+O(q))
   +2 \zeta^{9 \nu} q^{\frac{5}{26}} (1+O(q))+\\
  &+2 \zeta^{12 \nu} q^{\frac{11}{26}} (-1+O(q))+2 \zeta^{10 \nu} q^{\frac{7}{26}} (-1+O(q))
   +2 \zeta^{4 \nu} q^{\frac{47}{26}} (-1+O(q))+\\
  &+\zeta^{2 \nu} q^{\frac{17}{26}} (6+O(q))+\zeta^{5 \nu} q^{\frac{23}{26}} (8+O(q))
   +\zeta^{6 \nu} q^{\frac{25}{26}} (-1+O(q))+\\
  &+\zeta^{11 \nu} q^{\frac{9}{26}} (-3+O(q))+\zeta^{8 \nu} q^{\frac{29}{26}} (2+O(q))
   +\zeta^{7 \nu} q^{\frac{1}{26}} (1+O(q))]^8.
\endaligned$$
We will calculate the $q$-term which is the lowest degree.
(1) For the partition $26=7 \cdot 1+19$, the corresponding term is
$$\left(\matrix 8\\ 7, 1 \endmatrix\right) (\zeta^{7 \nu} q^{\frac{1}{26}})^7 \cdot
  2 \zeta^{3 \nu} q^{\frac{19}{26}} \cdot (-3)=-48 q.$$
(2) For the partition $26=6 \cdot 1+9+11$, the corresponding term is
$$\left(\matrix 8\\ 6, 1, 1 \endmatrix\right) (\zeta^{7 \nu} q^{\frac{1}{26}})^6 \cdot
  (\zeta^{11 \nu} q^{\frac{9}{26}}) \cdot (-3) \cdot 2 \zeta^{12 \nu} q^{\frac{11}{26}} (-1)=336 q.$$
(3) For the partition $26=6 \cdot 1+7+13$, the corresponding term is
$$\left(\matrix 8\\ 6, 1, 1 \endmatrix\right) (\zeta^{7 \nu} q^{\frac{1}{26}})^6 \cdot
  2 \zeta^{10 \nu} q^{\frac{7}{26}} \cdot (-1) \cdot q^{\frac{1}{2}} (-1)=112 q.$$
(4) For the partition $26=5 \cdot 1+5+7+9$, the corresponding term is
$$\left(\matrix 8\\ 5, 1, 1, 1 \endmatrix\right) (\zeta^{7 \nu} q^{\frac{1}{26}})^5 \cdot
  2 \zeta^{9 \nu} q^{\frac{5}{26}} \cdot 2 \zeta^{10 \nu} q^{\frac{7}{26}} (-1)
  \cdot \zeta^{11 \nu} q^{\frac{9}{26}} \cdot (-3)=4032 q.$$
(5) For the partition $26=5 \cdot 1+2 \cdot 5+11$, the corresponding term is
$$\left(\matrix 8\\ 5, 2, 1 \endmatrix\right) (\zeta^{7 \nu} q^{\frac{1}{26}})^5 \cdot
  (2 \zeta^{9 \nu} q^{\frac{5}{26}})^2 \cdot 2 \zeta^{12 \nu} q^{\frac{11}{26}} (-1)=-1344 q.$$
(6) For the partition $26=5 \cdot 1+3 \cdot 7$, the corresponding term is
$$\left(\matrix 8\\ 5, 3 \endmatrix\right) (\zeta^{7 \nu} q^{\frac{1}{26}})^5 \cdot
  (2 \zeta^{10 \nu} q^{\frac{7}{26}})^3 \cdot (-1)^3=-448 q.$$
(7) For the partition $26=4 \cdot 1+3 \cdot 5+7$, the corresponding term is
$$\left(\matrix 8\\ 4, 3, 1 \endmatrix\right) (\zeta^{7 \nu} q^{\frac{1}{26}})^4 \cdot
  (2 \zeta^{9 \nu} q^{\frac{5}{26}})^3 \cdot 2 \zeta^{10 \nu} q^{\frac{7}{26}} (-1)=-4480 q.$$
Hence, for $\Phi_{32}(y_1(z), \cdots, y_6(z))$ which is a modular form for $\Gamma(1)$ with weight $64$,
the lowest degree term is given by
$$(-48+336+112+4032-1344-448-4480) q \cdot q^{\frac{3}{24} \cdot 32}=-1840 q^5.$$
Thus,
$$\Phi_{32}(y_1(z), \cdots, y_6(z))=q^5 (-13 \cdot 1840+O(q)).$$
The leading term of $\Phi_{32}(y_1(z), \cdots, y_6(z))$ together with its weight $64$ suffice to identify
this modular form with $\Phi_{32}(y_1(z), \cdots, y_6(z))=-13 \cdot 1840 \Delta(z)^5 E_4(z)$.
Consequently,
$$\Phi_{32}(x_1(z), \cdots, x_6(z))=-13 \cdot 1840 \Delta(z)^5 E_4(z)/\eta(z)^{64}=-13 \cdot 1840 \eta(z)^8 \Delta(z)^2 E_4(z).$$

  Let
$$\Phi_{44}=w_0^{11}+w_1^{11}+\cdots+w_{12}^{11}+w_{\infty}^{11}.$$
As a polynomial in six variables, $\Phi_{44}(z_1, z_2, z_3, z_4, z_5, z_6)$ is a $G$-invariant polynomial.
Similarly as above, we can show that $\Phi_{44}(y_1(z), \cdots, y_6(z))$ is a modular form of weight $88$
for the full modular group $\Gamma(1)$. Moreover, we will show that it is a cusp form. In fact,
$$\aligned
  &\Phi_{44}(a_1(z), \cdots, a_6(z))=13^{11} q^{\frac{11}{2}} (1+O(q))^{11}+\\
  &+\sum_{\nu=0}^{12} [q^{\frac{1}{2}} (-1+O(q))+\\
  &+2 \zeta^{\nu} q^{\frac{41}{26}} (-3+O(q))+2 \zeta^{3 \nu} q^{\frac{19}{26}} (-3+O(q))
   +2 \zeta^{9 \nu} q^{\frac{5}{26}} (1+O(q))+\\
  &+2 \zeta^{12 \nu} q^{\frac{11}{26}} (-1+O(q))+2 \zeta^{10 \nu} q^{\frac{7}{26}} (-1+O(q))
   +2 \zeta^{4 \nu} q^{\frac{47}{26}} (-1+O(q))+\\
  &+\zeta^{2 \nu} q^{\frac{17}{26}} (6+O(q))+\zeta^{5 \nu} q^{\frac{23}{26}} (8+O(q))
   +\zeta^{6 \nu} q^{\frac{25}{26}} (-1+O(q))+\\
  &+\zeta^{11 \nu} q^{\frac{9}{26}} (-3+O(q))+\zeta^{8 \nu} q^{\frac{29}{26}} (2+O(q))
   +\zeta^{7 \nu} q^{\frac{1}{26}} (1+O(q))]^{11}.
\endaligned$$
We will calculate the $q^{\frac{3}{2}}$-term which is the lowest degree.
(1) For the partition $39=10 \cdot 1+29$, the corresponding term is
$$\left(\matrix 11\\ 10, 1 \endmatrix\right) (\zeta^{7 \nu} q^{\frac{1}{26}})^{10} \cdot
  \zeta^{8 \nu} q^{\frac{29}{26}} \cdot 2=22 q^{\frac{3}{2}}.$$
(2) For the partition $39=9 \cdot 1+13+17$, the corresponding term is
$$\left(\matrix 11\\ 9, 1, 1 \endmatrix\right) (\zeta^{7 \nu} q^{\frac{1}{26}})^{9} \cdot
  q^{\frac{1}{2}} \cdot (-1) \cdot \zeta^{2 \nu} q^{\frac{17}{26}} \cdot 6=-660 q^{\frac{3}{2}}.$$
(3) For the partition $39=9 \cdot 1+11+19$, the corresponding term is
$$\left(\matrix 11\\ 9, 1, 1 \endmatrix\right) (\zeta^{7 \nu} q^{\frac{1}{26}})^{9} \cdot
  2 \zeta^{3 \nu} q^{\frac{19}{26}} \cdot (-3) \cdot 2 \zeta^{12 \nu} q^{\frac{11}{26}} \cdot (-1)
  =1320 q^{\frac{3}{2}}.$$
(4) For the partition $39=9 \cdot 1+7+23$, the corresponding term is
$$\left(\matrix 11\\ 9, 1, 1 \endmatrix\right) (\zeta^{7 \nu} q^{\frac{1}{26}})^{9} \cdot
  2 \zeta^{10 \nu} q^{\frac{7}{26}} \cdot (-1) \cdot \zeta^{5 \nu} q^{\frac{23}{26}} \cdot 8
  =-1760 q^{\frac{3}{2}}.$$
(5) For the partition $39=9 \cdot 1+5+25$, the corresponding term is
$$\left(\matrix 11\\ 9, 1, 1 \endmatrix\right) (\zeta^{7 \nu} q^{\frac{1}{26}})^{9} \cdot
  2 \zeta^{9 \nu} q^{\frac{5}{26}} \cdot \zeta^{6 \nu} q^{\frac{25}{26}} \cdot (-1)=-220 q^{\frac{3}{2}}.$$
(6) For the partition $39=8 \cdot 1+13+2 \cdot 9$, the corresponding term is
$$\left(\matrix 11\\ 8, 1, 2 \endmatrix\right) (\zeta^{7 \nu} q^{\frac{1}{26}})^{8} \cdot
  q^{\frac{1}{2}} \cdot (-1) \cdot (\zeta^{11 \nu} q^{\frac{9}{26}})^2 \cdot (-3)^2=-4455 q^{\frac{3}{2}}.$$
(7) For the partition $39=8 \cdot 1+13+7+11$, the corresponding term is
$$\left(\matrix 11\\ 8, 1, 1, 1 \endmatrix\right) (\zeta^{7 \nu} q^{\frac{1}{26}})^{8} \cdot
  q^{\frac{1}{2}} \cdot (-1) \cdot 2 \zeta^{10 \nu} q^{\frac{7}{26}} \cdot (-1) \cdot 2
  \zeta^{12 \nu} q^{\frac{11}{26}} \cdot (-1)=-3960 q^{\frac{3}{2}}.$$
(8) For the partition $39=8 \cdot 1+2 \cdot 13+5$, the corresponding term is
$$\left(\matrix 11\\ 8, 2, 1 \endmatrix\right) (\zeta^{7 \nu} q^{\frac{1}{26}})^{8} \cdot
  (q^{\frac{1}{2}} \cdot (-1))^2 \cdot 2 \zeta^{9 \nu} q^{\frac{5}{26}}=990 q^{\frac{3}{2}}.$$
(9) For the partition $39=8 \cdot 1+19+5+7$, the corresponding term is
$$\left(\matrix 11\\ 8, 1, 1, 1 \endmatrix\right) (\zeta^{7 \nu} q^{\frac{1}{26}})^{8} \cdot
  2 \zeta^{3 \nu} q^{\frac{19}{26}} \cdot (-3) \cdot 2 \zeta^{9 \nu} q^{\frac{5}{26}} \cdot
  2 \zeta^{10 \nu} q^{\frac{7}{26}} \cdot (-1)=23760 q^{\frac{3}{2}}.$$
(10) For the partition $39=8 \cdot 1+2 \cdot 11+9$, the corresponding term is
$$\left(\matrix 11\\ 8, 2, 1 \endmatrix\right) (\zeta^{7 \nu} q^{\frac{1}{26}})^{8} \cdot
  (2 \zeta^{12 \nu} q^{\frac{11}{26}} \cdot (-1))^2 \cdot \zeta^{11 \nu} q^{\frac{9}{26}} \cdot (-3)=-5940 q^{\frac{3}{2}}.$$
(11) For the partition $39=8 \cdot 1+5+17+9$, the corresponding term is
$$\left(\matrix 11\\ 8, 1, 1, 1 \endmatrix\right) (\zeta^{7 \nu} q^{\frac{1}{26}})^{8} \cdot
  2 \zeta^{9 \nu} q^{\frac{5}{26}} \cdot \zeta^{2 \nu} q^{\frac{17}{26}} \cdot 6 \cdot \zeta^{11 \nu}
  q^{\frac{9}{26}} \cdot (-3)=-35640 q^{\frac{3}{2}}.$$
(12) For the partition $39=8 \cdot 1+2 \cdot 7+17$, the corresponding term is
$$\left(\matrix 11\\ 8, 2, 1 \endmatrix\right) (\zeta^{7 \nu} q^{\frac{1}{26}})^{8} \cdot
  (2 \zeta^{10 \nu} q^{\frac{7}{26}} \cdot (-1))^2 \cdot \zeta^{2 \nu} q^{\frac{17}{26}}
  \cdot 6=11880 q^{\frac{3}{2}}.$$
(13) For the partition $39=7 \cdot 1+13+5+2 \cdot 7$, the corresponding term is
$$\left(\matrix 11\\ 7, 1, 1, 2 \endmatrix\right) (\zeta^{7 \nu} q^{\frac{1}{26}})^{7} \cdot
  q^{\frac{1}{2}} \cdot (-1) \cdot 2 \zeta^{9 \nu} q^{\frac{5}{26}} \cdot (2 \zeta^{10 \nu}
  q^{\frac{7}{26}} \cdot (-1))^2=-31680 q^{\frac{3}{2}}.$$
(14) For the partition $39=7 \cdot 1+13+2 \cdot 5+9$, the corresponding term is
$$\left(\matrix 11\\ 7, 1, 2, 1 \endmatrix\right) (\zeta^{7 \nu} q^{\frac{1}{26}})^{7} \cdot
  q^{\frac{1}{2}} \cdot (-1) \cdot (2 \zeta^{9 \nu} q^{\frac{5}{26}})^2 \cdot \zeta^{11 \nu}
  q^{\frac{9}{26}} \cdot (-3)=47520 q^{\frac{3}{2}}.$$
(15) For the partition $39=7 \cdot 1+5+11+7+9$, the corresponding term is
$$\aligned
 &\left(\matrix 11\\ 7, 1, 1, 1, 1 \endmatrix\right) (\zeta^{7 \nu} q^{\frac{1}{26}})^{7} \cdot
  2 \zeta^{9 \nu} q^{\frac{5}{26}} \cdot 2 \zeta^{12 \nu} q^{\frac{11}{26}} \cdot (-1) \cdot 2
  \zeta^{10 \nu} q^{\frac{7}{26}} \cdot (-1) \cdot \zeta^{11 \nu} q^{\frac{9}{26}} \cdot (-3)\\
=&-190080 q^{\frac{3}{2}}.
\endaligned$$
(16) For the partition $39=7 \cdot 1+17+3 \cdot 5$, the corresponding term is
$$\left(\matrix 11\\ 7, 3, 1 \endmatrix\right) (\zeta^{7 \nu} q^{\frac{1}{26}})^{7} \cdot
  \zeta^{2 \nu} q^{\frac{17}{26}} \cdot 6 \cdot (2 \zeta^{9 \nu} q^{\frac{5}{26}})^3=63360 q^{\frac{3}{2}}.$$
(17) For the partition $39=7 \cdot 1+5+3 \cdot 9$, the corresponding term is
$$\left(\matrix 11\\ 7, 3, 1 \endmatrix\right) (\zeta^{7 \nu} q^{\frac{1}{26}})^{7} \cdot
  2 \zeta^{9 \nu} q^{\frac{5}{26}} \cdot (\zeta^{11 \nu} q^{\frac{9}{26}} \cdot (-3))^3=-71280 q^{\frac{3}{2}}.$$
(18) For the partition $39=7 \cdot 1+2 \cdot 5+2 \cdot 11$, the corresponding term is
$$\left(\matrix 11\\ 7, 2, 2 \endmatrix\right) (\zeta^{7 \nu} q^{\frac{1}{26}})^{7} \cdot
  (2 \zeta^{9 \nu} q^{\frac{5}{26}})^2 \cdot (2 \zeta^{12 \nu} q^{\frac{11}{26}} \cdot (-1))^2=31680 q^{\frac{3}{2}}.$$
(19) For the partition $39=7 \cdot 1+11+3 \cdot 7$, the corresponding term is
$$\left(\matrix 11\\ 7, 3, 1 \endmatrix\right) (\zeta^{7 \nu} q^{\frac{1}{26}})^{7} \cdot
  2 \zeta^{12 \nu} q^{\frac{11}{26}} \cdot (-1) \cdot (2 \zeta^{10 \nu} q^{\frac{7}{26}} \cdot (-1))^3=21120 q^{\frac{3}{2}}.$$
(20) For the partition $39=7 \cdot 1+2 \cdot 9+2 \cdot 7$, the corresponding term is
$$\left(\matrix 11\\ 7, 2, 2 \endmatrix\right) (\zeta^{7 \nu} q^{\frac{1}{26}})^{7} \cdot
  (\zeta^{11 \nu} q^{\frac{9}{26}} \cdot (-3))^2 \cdot (2 \zeta^{10 \nu} q^{\frac{7}{26}} \cdot (-1))^2=71280 q^{\frac{3}{2}}.$$
(21) For the partition $39=6 \cdot 1+13+4 \cdot 5$, the corresponding term is
$$\left(\matrix 11\\ 6, 1, 4 \endmatrix\right) (\zeta^{7 \nu} q^{\frac{1}{26}})^{6} \cdot
  q^{\frac{1}{2}} \cdot (-1) \cdot (2 \zeta^{9 \nu} q^{\frac{5}{26}})^4=-36960 q^{\frac{3}{2}}.$$
(22) For the partition $39=6 \cdot 1+5+4 \cdot 7$, the corresponding term is
$$\left(\matrix 11\\ 6, 1, 4 \endmatrix\right) (\zeta^{7 \nu} q^{\frac{1}{26}})^{6} \cdot
  2 \zeta^{9 \nu} q^{\frac{5}{26}} \cdot (2 \zeta^{10 \nu} q^{\frac{7}{26}} \cdot (-1))^4=73920 q^{\frac{3}{2}}.$$
(23) For the partition $39=6 \cdot 1+2 \cdot 5+2 \cdot 7+9$, the corresponding term is
$$\left(\matrix 11\\ 6, 2, 2, 1 \endmatrix\right) (\zeta^{7 \nu} q^{\frac{1}{26}})^{6} \cdot
  (2 \zeta^{9 \nu} q^{\frac{5}{26}})^2 \cdot (2 \zeta^{10 \nu} q^{\frac{7}{26}} \cdot (-1))^2
  \cdot \zeta^{11 \nu} q^{\frac{9}{26}} \cdot (-3)=-665280 q^{\frac{3}{2}}.$$
(24) For the partition $39=6 \cdot 1+3 \cdot 5+2 \cdot 9$, the corresponding term is
$$\left(\matrix 11\\ 6, 3, 2 \endmatrix\right) (\zeta^{7 \nu} q^{\frac{1}{26}})^{6} \cdot
  (2 \zeta^{9 \nu} q^{\frac{5}{26}})^3 \cdot (\zeta^{11 \nu} q^{\frac{9}{26}} \cdot (-3))^2=332640 q^{\frac{3}{2}}.$$
(25) For the partition $39=6 \cdot 1+3 \cdot 5+7+11$, the corresponding term is
$$\left(\matrix 11\\ 6, 3, 1, 1 \endmatrix\right) (\zeta^{7 \nu} q^{\frac{1}{26}})^{6} \cdot
  (2 \zeta^{9 \nu} q^{\frac{5}{26}})^3 \cdot 2 \zeta^{10 \nu} q^{\frac{7}{26}} \cdot (-1) \cdot
  2 \zeta^{12 \nu} q^{\frac{11}{26}} \cdot (-1)=295680 q^{\frac{3}{2}}.$$
(26) For the partition $39=5 \cdot 1+5 \cdot 5+9$, the corresponding term is
$$\left(\matrix 11\\ 5, 5, 1 \endmatrix\right) (\zeta^{7 \nu} q^{\frac{1}{26}})^{5} \cdot
  (2 \zeta^{9 \nu} q^{\frac{5}{26}})^5 \cdot \zeta^{11 \nu} q^{\frac{9}{26}} \cdot (-3)=-266112 q^{\frac{3}{2}}.$$
(27) For the partition $39=5 \cdot 1+4 \cdot 5+2 \cdot 7$, the corresponding term is
$$\left(\matrix 11\\ 5, 4, 2 \endmatrix\right) (\zeta^{7 \nu} q^{\frac{1}{26}})^{5} \cdot
  (2 \zeta^{9 \nu} q^{\frac{5}{26}})^4 \cdot (2 \zeta^{10 \nu} q^{\frac{7}{26}} \cdot (-1))^2=443520 q^{\frac{3}{2}}.$$
(28) For the partition $39=4 \cdot 1+7 \cdot 5$, the corresponding term is
$$\left(\matrix 11\\ 4, 7 \endmatrix\right) (\zeta^{7 \nu} q^{\frac{1}{26}})^{4} \cdot
  (2 \zeta^{9 \nu} q^{\frac{5}{26}})^7=42240 q^{\frac{3}{2}}.$$
Hence, for $\Phi_{44}(y_1(z), \cdots, y_6(z))$ which is a modular form for $\Gamma(1)$ with weight $88$,
the lowest degree term is given by
$$\aligned
 &(22-660+1320-1760-220-4455-3960+990+23760-5940-35640+\\
 &+11880-31680+47520-190080+63360-71280+31680+21120+71280+\\
 &-36960+73920-665280+332640+295680-266112+443520+42240) q^{\frac{3}{2}} \cdot q^{\frac{3}{24} \cdot 44}\\
=&146905 q^7.
\endaligned$$
Thus,
$$\Phi_{44}(y_1(z), \cdots, y_6(z))=q^7 (13 \cdot 146905+O(q)).$$
The leading term of $\Phi_{44}(y_1(z), \cdots, y_6(z))$ together with its weight $88$ suffice to identify
this modular form with $\Phi_{44}(y_1(z), \cdots, y_6(z))=13 \cdot 146905 \Delta(z)^7 E_4(z)$.
Consequently,
$$\Phi_{44}(x_1(z), \cdots, x_6(z))=13 \cdot 146905 \Delta(z)^7 E_4(z)/\eta(z)^{88}=13 \cdot 146905 \eta(z)^8 \Delta(z)^3 E_4(z).$$

  Let
$$\Phi_{42}=\delta_0^7+\delta_1^7+\cdots+\delta_{12}^7+\delta_{\infty}^7.$$
As a polynomial in six variables, $\Phi_{42}(z_1, z_2, z_3, z_4, z_5, z_6)$ is
a $G$-invariant polynomial. Similarly as above, we can show that $\Phi_{42}(x_1(z), \cdots, x_6(z))$ is a
modular form of weight $42$ for the full modular group $\Gamma(1)$. Moreover, we will show that it is a
cusp form. In fact,
$$\aligned
  &\Phi_{42}(a_1(z), \cdots, a_6(z))=13^{14} q^{\frac{49}{4}} (1+O(q))^7+\\
  &+\sum_{\nu=0}^{12} [-13 q^{\frac{7}{4}} (1+O(q))+\\
  &+\zeta^{\nu} q^{\frac{43}{52}} (13+O(q))+\zeta^{2 \nu} q^{\frac{47}{52}} (-22+O(q))
   +\zeta^{3 \nu} q^{\frac{51}{52}} (-21+O(q))+\\
  &+\zeta^{4 \nu} q^{\frac{3}{52}} (-1+O(q))+\zeta^{5 \nu} q^{\frac{7}{52}} (2+O(q))
   +\zeta^{6 \nu} q^{\frac{11}{52}} (2+O(q))+\\
  &+\zeta^{7 \nu} q^{\frac{15}{52}} (-2+O(q))+\zeta^{8 \nu} q^{\frac{19}{52}} (-8+O(q))
   +\zeta^{9 \nu} q^{\frac{23}{52}} (6+O(q))+\\
  &+\zeta^{10 \nu} q^{\frac{27}{52}} (1+O(q))+\zeta^{11 \nu} q^{\frac{31}{52}} (-8+O(q))
   +\zeta^{12 \nu} q^{\frac{35}{52}} (17+O(q))]^7.
\endaligned$$
Note that $3, 7, 11, 15, 19, 23, 27, 31, 35 \equiv 3$ (mod $4$), and $7(4k+3) \equiv 1$ (mod $4$),
but $39 \equiv 3$ (mod $4$). This implies that $q^{\frac{5}{4}}$-term is the lowest degree.
(1) For the partition $65=6 \cdot 3+47$, the corresponding term is
$$\left(\matrix 7\\ 6, 1 \endmatrix\right) (\zeta^{4 \nu} q^{\frac{3}{52}} \cdot (-1))^6 \cdot
  \zeta^{2 \nu} q^{\frac{47}{52}} \cdot (-22)=-154 q^{\frac{5}{4}}.$$
(2) For the partition $65=5 \cdot 3+7+43$, the corresponding term is
$$\left(\matrix 7\\ 5, 1, 1 \endmatrix\right) (\zeta^{4 \nu} q^{\frac{3}{52}} \cdot (-1))^5 \cdot
  \zeta^{5 \nu} q^{\frac{7}{52}} \cdot 2 \cdot \zeta^{\nu} q^{\frac{43}{52}} \cdot 13=-1092 q^{\frac{5}{4}}.$$
(3) For the partition $65=5 \cdot 3+15+35$, the corresponding term is
$$\left(\matrix 7\\ 5, 1, 1 \endmatrix\right) (\zeta^{4 \nu} q^{\frac{3}{52}} \cdot (-1))^5 \cdot
  \zeta^{7 \nu} q^{\frac{15}{52}} \cdot (-2) \cdot \zeta^{12 \nu} q^{\frac{35}{52}} \cdot 17=1428 q^{\frac{5}{4}}.$$
(4) For the partition $65=5 \cdot 3+19+31$, the corresponding term is
$$\left(\matrix 7\\ 5, 1, 1 \endmatrix\right) (\zeta^{4 \nu} q^{\frac{3}{52}} \cdot (-1))^5 \cdot
  \zeta^{8 \nu} q^{\frac{19}{52}} \cdot (-8) \cdot \zeta^{11 \nu} q^{\frac{31}{52}} \cdot (-8)=-2688 q^{\frac{5}{4}}.$$
(5) For the partition $65=5 \cdot 3+23+27$, the corresponding term is
$$\left(\matrix 7\\ 5, 1, 1 \endmatrix\right) (\zeta^{4 \nu} q^{\frac{3}{52}} \cdot (-1))^5 \cdot
  \zeta^{9 \nu} q^{\frac{23}{52}} \cdot 6 \cdot \zeta^{10 \nu} q^{\frac{27}{52}}=-252 q^{\frac{5}{4}}.$$
(6) For the partition $65=4 \cdot 3+7+11+35$, the corresponding term is
$$\left(\matrix 7\\ 4, 1, 1, 1 \endmatrix\right) (\zeta^{4 \nu} q^{\frac{3}{52}} \cdot (-1))^4 \cdot
  \zeta^{5 \nu} q^{\frac{7}{52}} \cdot 2 \cdot \zeta^{6 \nu} q^{\frac{11}{52}} \cdot 2 \cdot
  \zeta^{12 \nu} q^{\frac{35}{52}} \cdot 17=14280 q^{\frac{5}{4}}.$$
(7) For the partition $65=4 \cdot 3+7+15+31$, the corresponding term is
$$\left(\matrix 7\\ 4, 1, 1, 1 \endmatrix\right) (\zeta^{4 \nu} q^{\frac{3}{52}} \cdot (-1))^4 \cdot
  \zeta^{5 \nu} q^{\frac{7}{52}} \cdot 2 \cdot \zeta^{7 \nu} q^{\frac{15}{52}} \cdot (-2) \cdot
  \zeta^{11 \nu} q^{\frac{31}{52}} \cdot (-8)=6720 q^{\frac{5}{4}}.$$
(8) For the partition $65=4 \cdot 3+7+19+27$, the corresponding term is
$$\left(\matrix 7\\ 4, 1, 1, 1 \endmatrix\right) (\zeta^{4 \nu} q^{\frac{3}{52}} \cdot (-1))^4 \cdot
  \zeta^{5 \nu} q^{\frac{7}{52}} \cdot 2 \cdot \zeta^{8 \nu} q^{\frac{19}{52}} \cdot (-8) \cdot
  \zeta^{10 \nu} q^{\frac{27}{52}}=-3360 q^{\frac{5}{4}}.$$
(9) For the partition $65=4 \cdot 3+7+2 \cdot 23$, the corresponding term is
$$\left(\matrix 7\\ 4, 1, 2 \endmatrix\right) (\zeta^{4 \nu} q^{\frac{3}{52}} \cdot (-1))^4 \cdot
  \zeta^{5 \nu} q^{\frac{7}{52}} \cdot 2 \cdot (\zeta^{9 \nu} q^{\frac{23}{52}} \cdot 6)^2=7560 q^{\frac{5}{4}}.$$
(10) For the partition $65=4 \cdot 3+11+15+27$, the corresponding term is
$$\left(\matrix 7\\ 4, 1, 1, 1 \endmatrix\right) (\zeta^{4 \nu} q^{\frac{3}{52}} \cdot (-1))^4 \cdot
  \zeta^{6 \nu} q^{\frac{11}{52}} \cdot 2 \cdot \zeta^{7 \nu} q^{\frac{15}{52}} \cdot (-2) \cdot
  \zeta^{10 \nu} q^{\frac{27}{52}}=-840 q^{\frac{5}{4}}.$$
(11) For the partition $65=4 \cdot 3+11+19+23$, the corresponding term is
$$\left(\matrix 7\\ 4, 1, 1, 1 \endmatrix\right) (\zeta^{4 \nu} q^{\frac{3}{52}} \cdot (-1))^4 \cdot
  \zeta^{6 \nu} q^{\frac{11}{52}} \cdot 2 \cdot \zeta^{8 \nu} q^{\frac{19}{52}} \cdot (-8) \cdot
  \zeta^{9 \nu} q^{\frac{23}{52}} \cdot 6=-20160 q^{\frac{5}{4}}.$$
(12) For the partition $65=4 \cdot 3+2 \cdot 11+31$, the corresponding term is
$$\left(\matrix 7\\ 4, 2, 1 \endmatrix\right) (\zeta^{4 \nu} q^{\frac{3}{52}} \cdot (-1))^4 \cdot
  (\zeta^{6 \nu} q^{\frac{11}{52}} \cdot 2)^2 \cdot \zeta^{11 \nu} q^{\frac{31}{52}} \cdot (-8)=-3360 q^{\frac{5}{4}}.$$
(13) For the partition $65=4 \cdot 3+15+2 \cdot 19$, the corresponding term is
$$\left(\matrix 7\\ 4, 1, 2 \endmatrix\right) (\zeta^{4 \nu} q^{\frac{3}{52}} \cdot (-1))^4 \cdot
  \zeta^{7 \nu} q^{\frac{15}{52}} \cdot (-2) \cdot (\zeta^{8 \nu} q^{\frac{19}{52}} \cdot (-8))^2=-13440 q^{\frac{5}{4}}.$$
(14) For the partition $65=4 \cdot 3+2 \cdot 15+23$, the corresponding term is
$$\left(\matrix 7\\ 4, 2, 1 \endmatrix\right) (\zeta^{4 \nu} q^{\frac{3}{52}} \cdot (-1))^4 \cdot
  (\zeta^{7 \nu} q^{\frac{15}{52}} \cdot (-2))^2 \cdot \zeta^{9 \nu} q^{\frac{23}{52}} \cdot 6=2520 q^{\frac{5}{4}}.$$
(15) For the partition $65=3 \cdot 3+7+11+2 \cdot 19$, the corresponding term is
$$\left(\matrix 7\\ 3, 1, 1, 2 \endmatrix\right) (\zeta^{4 \nu} q^{\frac{3}{52}} \cdot (-1))^3 \cdot
  \zeta^{5 \nu} q^{\frac{7}{52}} \cdot 2 \cdot \zeta^{6 \nu} q^{\frac{11}{52}} \cdot 2 \cdot (\zeta^{8 \nu}
  q^{\frac{19}{52}} \cdot (-8))^2=-107520 q^{\frac{5}{4}}.$$
(16) For the partition $65=3 \cdot 3+7+11+15+23$, the corresponding term is
$$\left(\matrix 7\\ 3, 1, 1, 1, 1 \endmatrix\right) (\zeta^{4 \nu} q^{\frac{3}{52}} \cdot (-1))^3 \cdot
  \zeta^{5 \nu} q^{\frac{7}{52}} \cdot 2 \cdot \zeta^{6 \nu} q^{\frac{11}{52}} \cdot 2 \cdot \zeta^{7 \nu}
  q^{\frac{15}{52}} \cdot (-2) \cdot \zeta^{9 \nu} q^{\frac{23}{52}} \cdot 6=40320 q^{\frac{5}{4}}.$$
(17) For the partition $65=3 \cdot 3+7+2 \cdot 11+27$, the corresponding term is
$$\left(\matrix 7\\ 3, 1, 2, 1 \endmatrix\right) (\zeta^{4 \nu} q^{\frac{3}{52}} \cdot (-1))^3 \cdot
  \zeta^{5 \nu} q^{\frac{7}{52}} \cdot 2 \cdot (\zeta^{6 \nu} q^{\frac{11}{52}} \cdot 2)^2 \cdot
  \zeta^{10 \nu} q^{\frac{27}{52}}=-3360 q^{\frac{5}{4}}.$$
(18) For the partition $65=3 \cdot 3+7+2 \cdot 15+19$, the corresponding term is
$$\left(\matrix 7\\ 3, 1, 2, 1 \endmatrix\right) (\zeta^{4 \nu} q^{\frac{3}{52}} \cdot (-1))^3 \cdot
  \zeta^{5 \nu} q^{\frac{7}{52}} \cdot 2 \cdot (\zeta^{7 \nu} q^{\frac{15}{52}} \cdot (-2))^2 \cdot
  \zeta^{8 \nu} q^{\frac{19}{52}} \cdot(-8)=26880 q^{\frac{5}{4}}.$$
(19) For the partition $65=3 \cdot 3+2 \cdot 7+11+31$, the corresponding term is
$$\left(\matrix 7\\ 3, 2, 1, 1 \endmatrix\right) (\zeta^{4 \nu} q^{\frac{3}{52}} \cdot (-1))^3 \cdot
  (\zeta^{5 \nu} q^{\frac{7}{52}} \cdot 2)^2 \cdot \zeta^{6 \nu} q^{\frac{11}{52}} \cdot 2 \cdot
  \zeta^{11 \nu} q^{\frac{31}{52}} \cdot (-8)=26880 q^{\frac{5}{4}}.$$
(20) For the partition $65=3 \cdot 3+2 \cdot 7+15+27$, the corresponding term is
$$\left(\matrix 7\\ 3, 2, 1, 1 \endmatrix\right) (\zeta^{4 \nu} q^{\frac{3}{52}} \cdot (-1))^3 \cdot
  (\zeta^{5 \nu} q^{\frac{7}{52}} \cdot 2)^2 \cdot \zeta^{7 \nu} q^{\frac{15}{52}} \cdot (-2) \cdot
  \zeta^{10 \nu} q^{\frac{27}{52}}=3360 q^{\frac{5}{4}}.$$
(21) For the partition $65=3 \cdot 3+2 \cdot 7+19+23$, the corresponding term is
$$\left(\matrix 7\\ 3, 2, 1, 1 \endmatrix\right) (\zeta^{4 \nu} q^{\frac{3}{52}} \cdot (-1))^3 \cdot
  (\zeta^{5 \nu} q^{\frac{7}{52}} \cdot 2)^2 \cdot \zeta^{8 \nu} q^{\frac{19}{52}} \cdot (-8) \cdot
  \zeta^{9 \nu} q^{\frac{23}{52}} \cdot 6=80640 q^{\frac{5}{4}}.$$
(22) For the partition $65=3 \cdot 3+3 \cdot 7+35$, the corresponding term is
$$\left(\matrix 7\\ 3, 3, 1 \endmatrix\right) (\zeta^{4 \nu} q^{\frac{3}{52}} \cdot (-1))^3 \cdot
  (\zeta^{5 \nu} q^{\frac{7}{52}} \cdot 2)^3 \cdot \zeta^{12 \nu} q^{\frac{35}{52}} \cdot 17=-19040 q^{\frac{5}{4}}.$$
(23) For the partition $65=3 \cdot 3+11+3 \cdot 15$, the corresponding term is
$$\left(\matrix 7\\ 3, 1, 3 \endmatrix\right) (\zeta^{4 \nu} q^{\frac{3}{52}} \cdot (-1))^3 \cdot
  \zeta^{6 \nu} q^{\frac{11}{52}} \cdot 2 \cdot (\zeta^{7 \nu} q^{\frac{15}{52}} \cdot (-2))^3=2240 q^{\frac{5}{4}}.$$
(24) For the partition $65=3 \cdot 3+2 \cdot 11+15+19$, the corresponding term is
$$\left(\matrix 7\\ 3, 2, 1, 1 \endmatrix\right) (\zeta^{4 \nu} q^{\frac{3}{52}} \cdot (-1))^3 \cdot
  (\zeta^{6 \nu} q^{\frac{11}{52}} \cdot 2)^2 \cdot \zeta^{7 \nu} q^{\frac{15}{52}} \cdot (-2) \cdot
  \zeta^{8 \nu} q^{\frac{19}{52}} \cdot (-8)=-26880 q^{\frac{5}{4}}.$$
(25) For the partition $65=3 \cdot 3+3 \cdot 11+23$, the corresponding term is
$$\left(\matrix 7\\ 3, 3, 1 \endmatrix\right) (\zeta^{4 \nu} q^{\frac{3}{52}} \cdot (-1))^3 \cdot
  (\zeta^{6 \nu} q^{\frac{11}{52}} \cdot 2)^3 \cdot \zeta^{9 \nu} q^{\frac{23}{52}} \cdot 6=-6720 q^{\frac{5}{4}}.$$
(26) For the partition $65=2 \cdot 3+4 \cdot 7+31$, the corresponding term is
$$\left(\matrix 7\\ 2, 4, 1 \endmatrix\right) (\zeta^{4 \nu} q^{\frac{3}{52}} \cdot (-1))^2 \cdot
  (\zeta^{5 \nu} q^{\frac{7}{52}} \cdot 2)^4 \cdot \zeta^{11 \nu} q^{\frac{31}{52}} \cdot (-8)=-13440 q^{\frac{5}{4}}.$$
(27) For the partition $65=2 \cdot 3+3 \cdot 7+2 \cdot 19$, the corresponding term is
$$\left(\matrix 7\\ 2, 3, 2 \endmatrix\right) (\zeta^{4 \nu} q^{\frac{3}{52}} \cdot (-1))^2 \cdot
  (\zeta^{5 \nu} q^{\frac{7}{52}} \cdot 2)^3 \cdot (\zeta^{8 \nu} q^{\frac{19}{52}} \cdot (-8))^2=107520 q^{\frac{5}{4}}.$$
(28) For the partition $65=2 \cdot 3+3 \cdot 7+15+23$, the corresponding term is
$$\left(\matrix 7\\ 2, 3, 1, 1 \endmatrix\right) (\zeta^{4 \nu} q^{\frac{3}{52}} \cdot (-1))^2 \cdot
  (\zeta^{5 \nu} q^{\frac{7}{52}} \cdot 2)^3 \cdot \zeta^{7 \nu} q^{\frac{15}{52}} \cdot (-2) \cdot
  \zeta^{9 \nu} q^{\frac{23}{52}} \cdot 6=-40320 q^{\frac{5}{4}}.$$
(29) For the partition $65=2 \cdot 3+3 \cdot 7+11+27$, the corresponding term is
$$\left(\matrix 7\\ 2, 3, 1, 1 \endmatrix\right) (\zeta^{4 \nu} q^{\frac{3}{52}} \cdot (-1))^2 \cdot
  (\zeta^{5 \nu} q^{\frac{7}{52}} \cdot 2)^3 \cdot \zeta^{6 \nu} q^{\frac{11}{52}} \cdot 2 \cdot
  \zeta^{10 \nu} q^{\frac{27}{52}}=6720 q^{\frac{5}{4}}.$$
(30) For the partition $65=2 \cdot 3+2 \cdot 7+3 \cdot 15$, the corresponding term is
$$\left(\matrix 7\\ 2, 2, 3 \endmatrix\right) (\zeta^{4 \nu} q^{\frac{3}{52}} \cdot (-1))^2 \cdot
  (\zeta^{5 \nu} q^{\frac{7}{52}} \cdot 2)^2 \cdot (\zeta^{7 \nu} q^{\frac{15}{52}} \cdot (-2))^3=-6720 q^{\frac{5}{4}}.$$
(31) For the partition $65=2 \cdot 3+2 \cdot 7+11+15+19$, the corresponding term is
$$\aligned
 &\left(\matrix 7\\ 2, 2, 1, 1, 1 \endmatrix\right) (\zeta^{4 \nu} q^{\frac{3}{52}} \cdot (-1))^2 \cdot
  (\zeta^{5 \nu} q^{\frac{7}{52}} \cdot 2)^2 \cdot \zeta^{6 \nu} q^{\frac{11}{52}} \cdot 2 \cdot
  \zeta^{7 \nu} q^{\frac{15}{52}} \cdot (-2) \cdot \zeta^{8 \nu} q^{\frac{19}{52}} \cdot (-8)\\
=&161280 q^{\frac{5}{4}}.
\endaligned$$
(32) For the partition $65=2 \cdot 3+2 \cdot 7+2 \cdot 11+23$, the corresponding term is
$$\left(\matrix 7\\ 2, 2, 2, 1 \endmatrix\right) (\zeta^{4 \nu} q^{\frac{3}{52}} \cdot (-1))^2 \cdot
  (\zeta^{5 \nu} q^{\frac{7}{52}} \cdot 2)^2 \cdot (\zeta^{6 \nu} q^{\frac{11}{52}} \cdot 2)^2 \cdot
  \zeta^{9 \nu} q^{\frac{23}{52}} \cdot 6=60480 q^{\frac{5}{4}}.$$
(33) For the partition $65=2 \cdot 3+7+2 \cdot 11+2 \cdot 15$, the corresponding term is
$$\left(\matrix 7\\ 2, 1, 2, 2 \endmatrix\right) (\zeta^{4 \nu} q^{\frac{3}{52}} \cdot (-1))^2 \cdot
  \zeta^{5 \nu} q^{\frac{7}{52}} \cdot 2 \cdot (\zeta^{6 \nu} q^{\frac{11}{52}} \cdot 2)^2 \cdot
  (\zeta^{7 \nu} q^{\frac{15}{52}} \cdot (-2))^2=20160 q^{\frac{5}{4}}.$$
(34) For the partition $65=2 \cdot 3+7+3 \cdot 11+19$, the corresponding term is
$$\left(\matrix 7\\ 2, 1, 3, 1 \endmatrix\right) (\zeta^{4 \nu} q^{\frac{3}{52}} \cdot (-1))^2 \cdot
  \zeta^{5 \nu} q^{\frac{7}{52}} \cdot 2 \cdot (\zeta^{6 \nu} q^{\frac{11}{52}} \cdot 2)^3 \cdot
  \zeta^{8 \nu} q^{\frac{19}{52}} \cdot (-8)=-53760 q^{\frac{5}{4}}.$$
(35) For the partition $65=2 \cdot 3+4 \cdot 11+15$, the corresponding term is
$$\left(\matrix 7\\ 2, 4, 1 \endmatrix\right) (\zeta^{4 \nu} q^{\frac{3}{52}} \cdot (-1))^2 \cdot
  (\zeta^{6 \nu} q^{\frac{11}{52}} \cdot 2)^4 \cdot \zeta^{7 \nu} q^{\frac{15}{52}} \cdot (-2)=-3360 q^{\frac{5}{4}}.$$
(36) For the partition $65=3+5 \cdot 7+27$, the corresponding term is
$$\left(\matrix 7\\ 1, 5, 1 \endmatrix\right) \zeta^{4 \nu} q^{\frac{3}{52}} \cdot (-1) \cdot
  (\zeta^{5 \nu} q^{\frac{7}{52}} \cdot 2)^5 \cdot \zeta^{10 \nu} q^{\frac{27}{52}}=-1344 q^{\frac{5}{4}}.$$
(37) For the partition $65=3+4 \cdot 7+15+19$, the corresponding term is
$$\left(\matrix 7\\ 1, 4, 1, 1 \endmatrix\right) \zeta^{4 \nu} q^{\frac{3}{52}} \cdot (-1) \cdot
  (\zeta^{5 \nu} q^{\frac{7}{52}} \cdot 2)^4 \cdot \zeta^{7 \nu} q^{\frac{15}{52}} \cdot (-2) \cdot
  \zeta^{8 \nu} q^{\frac{19}{52}} \cdot (-8)=-53760 q^{\frac{5}{4}}.$$
(38) For the partition $65=3+4 \cdot 7+11+23$, the corresponding term is
$$\left(\matrix 7\\ 1, 4, 1, 1 \endmatrix\right) \zeta^{4 \nu} q^{\frac{3}{52}} \cdot (-1) \cdot
  (\zeta^{5 \nu} q^{\frac{7}{52}} \cdot 2)^4 \cdot \zeta^{6 \nu} q^{\frac{11}{52}} \cdot 2 \cdot
  \zeta^{9 \nu} q^{\frac{23}{52}} \cdot 6=-40320 q^{\frac{5}{4}}.$$
(39) For the partition $65=3+3 \cdot 7+11+2 \cdot 15$, the corresponding term is
$$\left(\matrix 7\\ 1, 3, 1, 2 \endmatrix\right) \zeta^{4 \nu} q^{\frac{3}{52}} \cdot (-1) \cdot
  (\zeta^{5 \nu} q^{\frac{7}{52}} \cdot 2)^3 \cdot \zeta^{6 \nu} q^{\frac{11}{52}} \cdot 2 \cdot
  (\zeta^{7 \nu} q^{\frac{15}{52}} \cdot (-2))^2=-26880 q^{\frac{5}{4}}.$$
(40) For the partition $65=3+3 \cdot 7+2 \cdot 11+19$, the corresponding term is
$$\left(\matrix 7\\ 1, 3, 2, 1 \endmatrix\right) \zeta^{4 \nu} q^{\frac{3}{52}} \cdot (-1) \cdot
  (\zeta^{5 \nu} q^{\frac{7}{52}} \cdot 2)^3 \cdot (\zeta^{6 \nu} q^{\frac{11}{52}} \cdot 2)^2 \cdot
  \zeta^{8 \nu} q^{\frac{19}{52}} \cdot (-8)=107520 q^{\frac{5}{4}}.$$
(41) For the partition $65=3+2 \cdot 7+3 \cdot 11+15$, the corresponding term is
$$\left(\matrix 7\\ 1, 2, 3, 1 \endmatrix\right) \zeta^{4 \nu} q^{\frac{3}{52}} \cdot (-1) \cdot
  (\zeta^{5 \nu} q^{\frac{7}{52}} \cdot 2)^2 \cdot (\zeta^{6 \nu} q^{\frac{11}{52}} \cdot 2)^3 \cdot
  \zeta^{7 \nu} q^{\frac{15}{52}} \cdot (-2)=26880 q^{\frac{5}{4}}.$$
(42) For the partition $65=3+7+5 \cdot 11$, the corresponding term is
$$\left(\matrix 7\\ 1, 1, 5 \endmatrix\right) \zeta^{4 \nu} q^{\frac{3}{52}} \cdot (-1) \cdot
  \zeta^{5 \nu} q^{\frac{7}{52}} \cdot 2 \cdot (\zeta^{6 \nu} q^{\frac{11}{52}} \cdot 2)^5=-2688 q^{\frac{5}{4}}.$$
(43) For the partition $65=6 \cdot 7+23$, the corresponding term is
$$\left(\matrix 7\\ 6, 1 \endmatrix\right) (\zeta^{5 \nu} q^{\frac{7}{52}} \cdot 2)^6 \cdot
  \zeta^{9 \nu} q^{\frac{23}{52}} \cdot 6=2688 q^{\frac{5}{4}}.$$
(44) For the partition $65=5 \cdot 7+2 \cdot 15$, the corresponding term is
$$\left(\matrix 7\\ 5, 2 \endmatrix\right) (\zeta^{5 \nu} q^{\frac{7}{52}} \cdot 2)^5 \cdot
  (\zeta^{7 \nu} q^{\frac{15}{52}} \cdot (-2))^2=2688 q^{\frac{5}{4}}.$$
(45) For the partition $65=5 \cdot 7+11+19$, the corresponding term is
$$\left(\matrix 7\\ 5, 1, 1 \endmatrix\right) (\zeta^{5 \nu} q^{\frac{7}{52}} \cdot 2)^5 \cdot
  \zeta^{6 \nu} q^{\frac{11}{52}} \cdot 2 \cdot \zeta^{8 \nu} q^{\frac{19}{52}} \cdot (-8)=-21504 q^{\frac{5}{4}}.$$
(46) For the partition $65=4 \cdot 7+2 \cdot 11+15$, the corresponding term is
$$\left(\matrix 7\\ 4, 2, 1 \endmatrix\right) (\zeta^{5 \nu} q^{\frac{7}{52}} \cdot 2)^4 \cdot
  (\zeta^{6 \nu} q^{\frac{11}{52}} \cdot 2)^2 \cdot \zeta^{7 \nu} q^{\frac{15}{52}} \cdot (-2)=-13440 q^{\frac{5}{4}}.$$
(47) For the partition $65=3 \cdot 7+4 \cdot 11$, the corresponding term is
$$\left(\matrix 7\\ 4, 3 \endmatrix\right) (\zeta^{5 \nu} q^{\frac{7}{52}} \cdot 2)^3 \cdot
  (\zeta^{6 \nu} q^{\frac{11}{52}} \cdot 2)^4=4480 q^{\frac{5}{4}}.$$
Hence, for $\Phi_{42}(x_1(z), \cdots, x_6(z))$ which is a modular form for $\Gamma(1)$ with weight $42$,
the lowest degree term is given by
$$\aligned
 &(-154-1092+1428-2688-252+14280+6720-3360+7560-840+\\
 &-20160-3360-13440+2520-107520+40320-3360+26880+26880+3360+\\
 &+80640-19040+2240-26880-6720-13440+107520-40320+6720-6720+\\
 &+161280+60480+20160-53760-3360-1344-53760-40320-26880+\\
 &+107520+26880-2688+2688+2688-21504-13440+4480) q^{\frac{5}{4}} \cdot q^{\frac{42}{24}}\\
 &=226842 q^3.
\endaligned$$
Thus,
$$\Phi_{42}(x_1(z), \cdots, x_6(z))=q^3 (13 \cdot 226842+O(q)).$$
The leading term of $\Phi_{42}(x_1(z), \cdots, x_6(z))$ together with its weight $42$ suffice to identify
this modular form with
$$\Phi_{42}(x_1(z), \cdots, x_6(z))=13 \cdot 226842 \Delta(z)^3 E_6(z).$$

  Up to a constant, we revise the definition of $\Phi_{32}$, $\Phi_{42}$ and $\Phi_{44}$:
$$\Phi_{32}=-\frac{1}{13 \cdot 1840} \left(\sum_{\nu=0}^{12} w_{\nu}^8+w_{\infty}^8\right), \quad
  \Phi_{42}=\frac{1}{13 \cdot 226842} \left(\sum_{\nu=0}^{12} \delta_{\nu}^7+\delta_{\infty}^7\right),\eqno{(4.1)}$$
$$\Phi_{44}=\frac{1}{13 \cdot 146905} \left(\sum_{\nu=0}^{12} w_{\nu}^{11}+w_{\infty}^{11}\right).\eqno{(4.2)}$$
Consequently,
$$\left\{\aligned
  \Phi_{32}(x_1(z), \cdots, x_6(z)) &=\eta(z)^8 \Delta(z)^2 E_4(z),\\
  \Phi_{42}(x_1(z), \cdots, x_6(z)) &=\Delta(z)^3 E_6(z),\\
  \Phi_{44}(x_1(z), \cdots, x_6(z)) &=\eta(z)^8 \Delta(z)^3 E_4(z).
\endaligned\right.\eqno{(4.3)}$$
From now on, we will use the following abbreviation $\Phi_j=\Phi_j(x_1(z), \cdots, x_6(z))$ for
$j=32, 42$ and $44$. The relations
$$j(z):=\frac{E_4(z)^3}{\Delta(z)}=\frac{\Phi_{32}^3}{\Phi_{12}^8}=\frac{\Phi_{44}^3}{\Phi_{12}^{11}}, \quad
  j(z)-1728=\frac{E_6(z)^2}{\Delta(z)}=\frac{\Phi_{42}^2}{\Phi_{12}^7}\eqno{(4.4)}$$
together with (3.20) give the equations
$$\left\{\aligned
  \Phi_{20}^3 \Phi_{12}^2-\Phi_{42}^2 &=1728 \Phi_{12}^7,\\
  \Phi_{32}^3-\Phi_{12}^5 \Phi_{18}^2 &=1728 \Phi_{12}^8,\\
  \Phi_{32}^3-\Phi_{12}^3 \Phi_{30}^2 &=1728 \Phi_{12}^8,\\
    \Phi_{32}^3-\Phi_{12} \Phi_{42}^2 &=1728 \Phi_{12}^8,\\
  \Phi_{44}^3-\Phi_{12}^8 \Phi_{18}^2 &=1728 \Phi_{12}^{11},\\
  \Phi_{44}^3-\Phi_{12}^6 \Phi_{30}^2 &=1728 \Phi_{12}^{11},\\
  \Phi_{44}^3-\Phi_{12}^4 \Phi_{42}^2 &=1728 \Phi_{12}^{11}.
\endaligned\right.\eqno{(4.5)}$$
This leads to the modular parametrizations of the following three exceptional singularities:
$$Q_{18}: x^3+y^8+y z^2=0, \quad E_{20}: x^3+y^{11}+z^2=0, \quad x^7+x^2 y^3+z^2=0.$$
Note that $Q_{18}$ and $E_{20}$ are two singularities in the pyramids of $14$ exceptional
singularities of modalities $2$ (see \cite{Ar4}, p.255). However, the last singularity does
not appear in the singularities with the number of moduli $m=0$, $1$ and $2$ (see\cite{Ar4}).
Moreover, the last three equations in (4.5) give a new analytic construction of solutions for the
Diophantine equation $x^p+y^q=z^r$ in the case $(p, q, r)=(2, 3, 11)$. Thus, we complete the proof of
Theorem 1.2 and Corollary 1.3.

\flushpar $\qquad \qquad \qquad \qquad \qquad \qquad \qquad \qquad
\qquad \qquad \qquad \qquad \qquad \qquad \qquad \qquad \qquad
\qquad \quad \boxed{}$

  It should be pointed out that in \cite{K}, \cite{K1}, \cite{K2}, \cite{K3}, \cite{KF1}, Klein connected the
simple groups $\text{PSL}(2, 5)$ and $\text{PSL}(2, 7)$ with singularities $E_8$ and $E_{12}$, respectively.
However, he could not connect the simple group $\text{PSL}(2, 11)$ with any singularities (see \cite{K4},
\cite{KF2}, p.413), especially the singularity $E_{20}$.  By Theorem 1.2, we connect the simple group
$\text{PSL}(2, 13)$ with singularity $E_{20}$ and solve this problem dating back to the works of Klein.
Moreover, by Theorem 1.1 and Theorem 1.2, we connect the simple group $\text{PSL}(2, 13)$ with many
singularities: $E_8$, $Q_{18}$, $E_{20}$ and $x^7+x^2 y^3+z^2=0$.

\vskip 2.0 cm

{\smc Department of Mathematics, Peking University}

{\smc Beijing 100871, P. R. China}

{\it E-mail address}: yanglei\@math.pku.edu.cn
\vskip 1.5 cm
\Refs

\item{[Ar1]} {\smc V. I. Arnold}, Critical points of smooth functions, in: {\it Proceedings of
             the International Congress of Mathematicians $($Vancouver, B. C., 1974$)$}, Vol. 1,
             19-39, Canad. Math. Congress, Montreal, Que., 1975.

\item{[Ar2]} {\smc V. I. Arnold}, Some open problems in the theory of singularities, translated from
              the Russian, in: {\it Singularities, Part 1 $($Arcata, Calif., 1981$)$}, 57-69, Proc. Sympos.
              Pure Math. {\bf 40}, Amer. Math. Soc., Providence, R.I., 1983.

\item{[Ar3]} {\smc V. I. Arnold}, Singularities of ray systems, in: {\it Proceedings of the International
             Congress of Mathematicians, Vol. 1, 2 $($Warsaw, 1983$)$}, 27-49, PWN, Warsaw, 1984.

\item{[Ar4]} {\smc V. I. Arnold, S. M. Gusein-Zade and  A. N. Varchenko}, {\it Singularities
             of Differentiable Maps}, Vol. I. {\it The Classification of Critical Points, Caustics
             and Wave Fronts}, Translated from the Russian by Ian Porteous and Mark Reynolds,
             Monographs in Mathematics, {\bf 82}, Birkh\"{a}user, 1985.

\item{[Ar5]} {\smc V. I. Arnold}, {\it Catastrophe Theory}, translated from the Russian by G. S. Wassermann,
             based on a translation by R. K. Thomas, Third edition, Springer-Verlag, Berlin, 1992.

\item{[Ar6]} {\smc V. I. Arnold}, {\it Arnold's Problems}, Springer-Verlag, Berlin; PHASIS, Moscow, 2004

\item{[At]} {\smc M. Atiyah}, The icosahedron, Math. Medley {\bf 18} (1990), 1-12.

\item{[B]} {\smc W. P. Barth, K. Hulek, C. A. M. Peters and A. Van de Ven}, {\it Compact complex surfaces},
              Second enlarged edition, Ergebnisse der Mathematik und ihrer Grenzgebiete, 3. Folge, Vol. {\bf 4},
              Springer-Verlag, Berlin, 2004.

\item{[Be]} {\smc F. Beukers}, The Diophantine equation $Ax^p+By^q=Cz^r$, Duke Math. J. {\bf 91} (1998), 61-88.

\item{[Br1]} {\smc E. Brieskorn}, Examples of singular normal complex spaces
             which are topological manifolds, Proc. Nat. Acad. Sci. U.S.A. {\bf 55}
             (1966), 1395-1397.

\item{[Br2]} {\smc E. Brieskorn}, Beispiele zur Differentialtopologie von Singularit\"{a}ten,
             Invent. Math. {\bf 2} (1966), 1-14.

\item{[Br3]} {\smc E. Brieskorn}, Singular elements of semi-simple algebraic groups, in: {\it Actes du
             Congr\`{e}s International des Math\'{e}maticiens (Nice, 1970)}, Tome 2, 279-284,
             Gauthier-Villars, Paris, 1971.

\item{[Br4]} {\smc E. Brieskorn}, The unfolding of exceptional singularities, in: {\it Leopoldina
             Symposium: Singularities (Th\"{u}ringen, 1978)}, Nova Acta Leopoldina (N.F.)
             {\bf 52} (1981), no. 240, 65-93.

\item{[Br5]} {\smc E. Brieskorn}, Singularities in the work of Friedrich Hirzebruch, in: {\it Surveys
              in differential geometry}, 17-60, Surv. Differ. Geom., {\bf 7}, Int. Press, Somerville, MA, 2000.

\item{[BPR]} {\smc E. Brieskorn, A. Pratoussevitch and F. Rothenh\"{a}usler}, The combinatorial
             geometry of singularities and Arnold's series $E$, $Z$, $Q$, Dedicated to Vladimir I.
             Arnold on the occasion of his 65th birthday, Mosc. Math. J. {\bf 3} (2003), 273-333, 741.

\item{[CIZ]} {\smc A. Cappelli, C. Itzykson and J.-B. Zuber}, The A-D-E classification of minimal and $A_1^{(1)}$
              conformal invariant theories, Comm. Math. Phys. {\bf 113} (1987), 1-26.

\item{[CC]} {\smc J. H. Conway, R. T. Curtis, S. P. Norton, R. A. Parker and R. A. Wilson}, {\it Atlas of
             Finite Groups, Maximal Subgroups and Ordinary Characters for Simple Groups}, Clarendon Press, Oxford, 1985.

\item{[CoM]} {\smc H. S. M. Coxeter and W.O. J. Moser}, {\it Generators and Relations for Discrete Groups},
             Third edition, Ergebnisse der Mathematik und ihrer Grenzgebiete, Band {\bf 14}, Springer-Verlag,
             New York-Heidelberg, 1972.

\item{[DG]} {\smc H. Darmon and A. Granville}, On the equations $z^m=F(x,y)$ and $Ax^p+By^q=Cz^r$,
             Bull. London Math. Soc. {\bf 27} (1995), 513-543.

\item{[D]} {\smc H. Darmon}, Faltings plus epsilon, Wiles plus epsilon, and the generalized Fermat equation,
           C. R. Math. Rep. Acad. Sci. Canada {\bf 19} (1997), 3-14, 64.

\item{[DM]} {\smc P. Doyle and C. McMullen}, Solving the quintic by iteration, Acta Math. {\bf 163} (1989), 151-180.

\item{[Du]} {\smc W. Duke}, Continued fractions and modular functions, Bull. Amer. Math. Soc. (N.S.) {\bf 42} (2005), 137-162.

\item{[Dur]} {\smc A. H. Durfee}, Fifteen characterizations of rational double points and simple critical points,
             Enseign. Math. (2) {\bf 25} (1979), 131-163.

\item{[DV]} {\smc P. Du Val}, On isolated singularities of surfaces which do not affect the conditions
            of adjunction, I, II, III, Math. Proc. Cambridge Philos. Soc. {\bf 30} (1934), 453-459,
            460-465, 483-491.

\item{[E]} {\smc J. Edwards}, A complete solution to $X^2+Y^3+Z^5=0$, J. Reine Angew. Math.
           {\bf 571} (2004), 213-236.

\item{[FK]} {\smc H. M. Farkas and I. Kra}, {\it Theta Constants,
            Riemann Surfaces and the Modular Group, An Introduction
            with Applications to Uniformization Theorems, Partition
            Identities and Combinatorial Number Theory}, Graduate
            Studies in Mathematics, {\bf 37}, American Mathematical
            Society, Providence, RI, 2001.

\item{[GV]} {\smc G. Gonzalez-Sprinberg and J.-L. Verdier}, Construction g\'{e}om\'{e}trique
            de la correspondance de McKay, Ann. Sci. \'{E}cole Norm. Sup. (4) {\bf 16} (1983),
            409-449.

\item{[Gr]} {\smc G.-M. Greuel}, Some aspects of Brieskorn's mathematical work,
            {\it Singularities $($Oberwolfach, 1996$)$}, xv-xxii, Progr. Math.,
            {\bf 162}, Birkh\"{a}user, Basel, 1998.

\item{[Hi1]} {\smc F. Hirzebruch}, The topology of normal singularities of an algebraic
            surface (d'apr\`{e}s un article de D. Mumford), S\'{e}minaire Bourbaki, 1962/63,
            Exp. 250, in: {\it Gesammelte Abhandlungen}, Bd. II, 1-7, Springer-Verlag, 1987.

\item{[Hi2]} {\smc F. Hirzebruch}, Singularities and exotic spheres, S\'{e}minaire Bourbaki,
            1966/67, Exp. 314, in: {\it Gesammelte Abhandlungen}, Bd. II, 70-80, Springer-Verlag, 1987.

\item{[Hi3]} {\smc F. Hirzebruch}, The icosahedron, in: {\it Gesammelte Abhandlungen}, Bd. II,
            656-661, Springer-Verlag, 1987.

\item{[Hi4]} {\smc F. Hirzebruch}, Zur Theorie der Mannigfaltigkeiten, in: {\it Gesammelte Abhandlungen}, Bd. I,
            673, Springer-Verlag, 1987.

\item{[HM]} {\smc F. Hirzebruch and K. Mayer}, {\it $O(n)$-Mannigfaltigkeiten, exotische Sph\"{a}ren und
            Singularit\"{a}ten}, Lecture Notes in Math., {\bf 57}, Springer-Verlag, 1968.

\item{[H]} {\smc N. Hitchin}, $E_6$, $E_7$, $E_8$, Clay academy lecture, 2005.

\item{[IN]} {\smc Y. Ito and I. Nakamura}, Hilbert schemes and simple singularities,
            in: {\it New trends in algebraic geometry $($Warwick, 1996$)$}, 151-233, London Math. Soc.
            Lecture Note Ser., {\bf 264}, Cambridge Univ. Press, Cambridge, 1999.

\item{[J1]} {\smc V. F. R. Jones}, Subfactors of type ${\text{II}}_1$ factors and related topics, in: {\it
            Proceedings of the International Congress of Mathematicians, Vol. 1, 2 $($Berkeley, Calif., 1986$)$},
            939-947, Amer. Math. Soc., Providence, RI, 1987.

\item{[J2]} {\smc V. F. R. Jones}, von Neumann algebras in mathematics and physics, in: {\it Proceedings of the
           International Congress of Mathematicians, Vol. I, II $($Kyoto, 1990$)$}, 121-138, Math. Soc. Japan, Tokyo, 1991.

\item{[J3]} {\smc V. F. R. Jones}, In and around the origin of quantum groups, in: {\it Prospects in mathematical physics},
            101-126, Contemp. Math. {\bf 437}, Amer. Math. Soc., Providence, RI, 2007.

\item{[KM]} {\smc M. Kervaire and J. Milnor}, Groups of homotopy spheres: I, Ann. of Math. (2)
            {\bf 77} (1963), 504-537.

\item{[KS]} {\smc R. C. Kirby and M. G. Scharlemann}, Eight faces of the Poincar\'{e} homology
             $3$-sphere, in: {\it Geometric topology $($Proc. Georgia Topology Conf., Athens, Ga.,
             1977$)$}, 113-146, Academic Press, New York-London, 1979.

\item{[K]} {\smc F. Klein}, {\it Lectures on the Icosahedron and
             the Solution of Equations of the Fifth Degree},
             Translated by G. G. Morrice, second and revised edition,
             Dover Publications, Inc., 1956.

\item{[K1]} {\smc F. Klein}, Ueber die Transformation der
            elliptischen Functionen und die Aufl\"{o}sung der
            Gleichungen f\"{u}nften Grades, Math. Ann. {\bf 14}
            (1879), 111-172, in: {\it Gesammelte Mathematische
            Abhandlungen}, Bd. III, 13-75, Springer-Verlag, Berlin,
            1923.

\item{[K2]} {\smc F. Klein}, Ueber die Transformation siebenter
            Ordnung der elliptischen Functionen, Math. Ann. {\bf
            14} (1879), 428-471, in: {\it Gesammelte Mathematische
            Abhandlungen}, Bd. III, 90-136, Springer-Verlag, Berlin,
            1923.

\item{[K3]} {\smc F. Klein}, Ueber die Aufl\"{o}sung gewisser
            Gleichungen vom siebenten und achten Grade, Math. Ann.
            {\bf 15} (1879), 251-282, in: {\it Gesammelte Mathematische
            Abhandlungen}, Bd. II, 390-438, Springer-Verlag, Berlin,
            1922.

\item{[K4]} {\smc F. Klein}, Ueber die Transformation elfter
            Ordnung der elliptischen Functionen, Math. Ann. {\bf
            15} (1879), 533-555, in: {\it Gesammelte Mathematische
            Abhandlungen}, Bd. III, 140-168, Springer-Verlag, Berlin,
            1923.

\item{[KF1]} {\smc F. Klein and R. Fricke}, {\it Vorlesungen \"{u}ber die Theorie der Elliptischen
            Modulfunctionen}, Vol. I, Leipzig, 1890.

\item{[KF2]} {\smc F. Klein and R. Fricke}, {\it Vorlesungen \"{u}ber die Theorie der Elliptischen
            Modulfunctionen}, Vol. II, Leipzig, 1892.

\item{[Ko1]} {\smc B. Kostant}, On finite subgroups of $\text{SU}(2)$, simple Lie algebras,
            and the McKay correspondence, Proc. Nat. Acad. Sci. U.S.A. {\bf 81} (1984), 5275-5277.

\item{[Ko2]} {\smc B. Kostant}, The McKay correspondence, the Coxeter element and representation theory,
            in: {\it The mathematical heritage of \'{E}lie Cartan $($Lyon, 1984$)$}, Ast\'{e}risque 1985,
            Num\'{e}ro Hors S\'{e}rie, 209-255.

\item{[Kr]}{\smc P. B. Kronheimer}, The construction of ALE spaces as hyper-K\"{a}hler quotients,
           J. Diff. Geom. {\bf 29} (1989), 665-683.

\item{[L]} {\smc J. H. Lindsey, II}, Finite linear groups of degree six, Canad. J. Math. {\bf 23} (1971), 771-790.

\item{[Mc1]} {\smc J. McKay}, Graphs, singularities, and finite groups, in: {\it The Santa
            Cruz Conference on Finite Groups $($Univ. California, Santa Cruz, Calif., 1979$)$},
            183-186, Proc. Sympos. Pure Math., {\bf 37}, Amer. Math. Soc., Providence, R.I., 1980.

\item{[Mc2]} {\smc J. McKay}, Cartan matrices, finite groups of quaternions, and Kleinian singularities,
             Proc. Amer. Math. Soc. {\bf 81} (1981), 153-154.

\item{[Mi1]} {\smc J. Milnor}, On manifolds homeomorphic to the $7$-sphere, Ann. of Math. (2) {\bf 64}
            (1956), 399-405.

\item{[Mi2]} {\smc J. Milnor}, Differentiable manifolds which are homotopy spheres, in: {\it Collected
             papers of John Milnor}, III. {\it Differential topology}, 65-88, American Mathematical
             Society, Providence, RI, 2007.

\item{[Mi3]} {\smc J. Milnor}, Differential topology, in: {\it Collected papers of John Milnor}, III.
            {\it Differential topology}, 123-141, American Mathematical Society, Providence, RI, 2007.

\item{[Mi4]} {\smc J. Milnor}, {\it Collected papers of John Milnor}, III. {\it Differential topology},
             American Mathematical Society, Providence, RI, 2007.

\item{[N]} {\smc H. Nakajima}, Geometric construction of representations of affine algebras, in:
           {\it Proceedings of the International Congress of Mathematicians, Vol. I $($Beijing, 2002$)$},
           423-438, Higher Ed. Press, Beijing, 2002.

\item{[Na]} {\smc I. Naruki}, $E_8$ und die bin\"{a}re Ikosaedergruppe, Invent. Math. {\bf 42} (1977), 273-283.

\item{[PSS]} {\smc B. Poonen, E. F. Schaefer and M. Stoll}, Twists of $X(7)$ and primitive solutions to
            $x^2+y^3=z^7$, Duke Math. J. {\bf 137} (2007), 103-158.

\item{[R]} {\smc M. Reid}, La correspondance de McKay, S\'{e}minaire Bourbaki, Vol. 1999/2000,
            Ast\'{e}risque No. {\bf 276} (2002), 53-72.

\item{[Sch]} {\smc H. A. Schwarz}, \"{U}ber diejenigen F\"{a}lle, in welchen die Gaussische hypergeometrische
            Reihe eine algebraische Function ihres vierten Elementes darstellt, J. Reine Angew. Math. {\bf 75} (1872),
            292-335.

\item{[Se]} {\smc J.-P. Serre}, Extensions icosa\'{e}driques, S\'{e}minaire de Th\'{e}orie des
            Nombres de Bordeaux 1979/80, No. 19, in: {\it {\OE}uvres}, Vol. III, 550-554, Springer-Verlag, 1986.

\item{[Sh]} {\smc G. Shimura}, Construction of class fields and zeta functions of algebraic curves, Ann. of Math.
            {\bf 85} (1967), 58-159, in: {\it Collected Papers}, Vol. II, 30-131, Springer-Verlag, 2002.

\item{[Sl1]} {\smc P. Slodowy}, {\it Simple singularities and simple algebraic groups},
             Lecture Notes in Math. {\bf 815}, Springer, Berlin, 1980.

\item{[Sl2]} {\smc P. Slodowy}, Platonic solids, Kleinian singularities, and Lie groups, in:
            {\it Algebraic geometry $($Ann Arbor, Mich., 1981$)$}, 102-138, Lecture Notes in Math.,
            {\bf 1008}, Springer, Berlin, 1983.

\item{[V]} {\smc D. Vogan}, The character table for $E_8$, http://www-math.mit.edu/$\sim$dav/e8wpiHOedit.pdf, 2011.

\item{[W]} {\smc E. Witten}, Singularities in string theory, in: {\it Proceedings of the International Congress of
            Mathematicians, Vol. I $($Beijing, 2002$)$}, 495-504, Higher Ed. Press, Beijing, 2002.

\item{[Y1]} {\smc L. Yang}, Exotic arithmetic structure on the first Hurwitz triplet, arXiv:1209.1783v5 [math.NT], 2013.

\item{[Y2]} {\smc L. Yang}, Dedekind $\eta$-function, Hauptmodul and invariant theory, arXiv:1407.3550v2 [math.NT], 2014.

\item{[Z]} {\smc J.-B. Zuber}, CFT, BCFT, ADE and all that, in: {\it Quantum symmetries in theoretical physics
            and mathematics $($Bariloche, 2000$)$}, 233-266, Contemp. Math., {\bf 294}, Amer. Math. Soc., Providence, RI, 2002.

\endRefs
\end{document}